\documentclass[journal]{IEEEtran}%

\usepackage{amsmath,amssymb,amsfonts}%
\usepackage{theorem}%
\usepackage{graphicx}%
\usepackage[table,dvipsnames]{xcolor}

\newtheorem{definition}{Definition}[section]%
\newtheorem{theorem}[definition]{Theorem}%
\newtheorem{proposition}[definition]{Proposition}%
\newtheorem{lemma}[definition]{Lemma}%
\newtheorem{assumption}[definition]{Assumption}%
\newtheorem{corollary}[definition]{Corollary}%
{\theorembodyfont{\rmfamily} \newtheorem{remark}[definition]{Remark}}%
{\theorembodyfont{\rmfamily} \newtheorem{example}[definition]{Example}}%

\newcommand{\tm}{\times}%
\newcommand{\trn}{^{\scriptscriptstyle \top}}%
\newcommand{\esssup}{\operatorname*{ess\;sup}}%

\newcommand{\id}{\mathrm{id}}%
\newcommand{\KL}{\mathcal{KL}}%

\newcommand{\ep}{\varepsilon}%
\newcommand{\C}{\mathbb{C}}%
\newcommand{\N}{\mathbb{N}}%
\newcommand{\R}{\mathbb{R}}%
\newcommand{\KC}{\mathcal{K}}%
\newcommand{\LC}{\mathcal{L}}%
\newcommand{\PC}{\mathcal{P}}%
\newcommand{\UC}{\mathcal{U}}%
\newcommand{\rmd}{\mathrm{d}}%
\newcommand{\rmD}{\mathrm{D}}%
\newcommand{\diag}{\mathrm{diag}}%
\newcommand{\inner}{\mathrm{int}}%
\newcommand{\rme}{\mathrm{e}}%
\newcommand{\normt}[1]{{\left\vert\kern-0.25ex\left\vert\kern-0.25ex\left\vert #1 
		\right\vert\kern-0.25ex\right\vert\kern-0.25ex\right\vert}}

\allowdisplaybreaks%

\begin{document}

\title{A Lyapunov-based small-gain theorem for infinite networks}%

\author{Christoph~Kawan, Andrii~Mironchenko, Abdalla~Swikir, Navid~Noroozi, Majid~Zamani,~\IEEEmembership{Senior Member,~IEEE}
	\thanks{C.~Kawan is with the Institute of Informatics, LMU Munich, Germany; e-mail: \texttt{christoph.kawan@lmu.de}. His work is supported by the German Research Foundation (DFG) through the grant ZA 873/4-1.}%
  \thanks{A.~Mironchenko is with Faculty of Computer Science and Mathematics, University of Passau, 94032 Passau, Germany; e-mail: \texttt{andrii.mironchenko@uni-passau.de}. His work is supported by the DFG through the grant MI 1886/2-1.}
	\thanks{N.~Noroozi is with the Institute of Informatics, LMU Munich, Germany; e-mail: \texttt{navid.noroozi@lmu.de}. His work is supported by the DFG through the grant WI 1458/16-1.}
\thanks{A.~Swikir is with Department of Electrical and Computer Engineering, Technical University of Munich, Germany; e-mail: \texttt{abdalla.swikir@tum.de}.}
	\thanks{M.~Zamani is with the Computer Science Department, University of Colorado Boulder, CO 80309, USA. 
M.~Zamani is also with the the Institute of Informatics, LMU Munich, Germany; email: {\tt\small majid.zamani@colorado.edu}. His work is supported in part by the DFG through the grant ZA 873/4-1 and the H2020 ERC Starting Grant AutoCPS (grant agreement No.~804639).}}

\maketitle%

\begin{abstract}
This paper presents a small-gain theorem for networks composed of a countably infinite number of finite-dimensional subsystems. Assuming that each subsystem is exponentially input-to-state stable, we show that if the gain operator, collecting all the information about the internal Lyapunov gains, has a spectral radius less than one, the overall infinite network is exponentially input-to-state stable. The effectiveness of our result is illustrated through several examples including nonlinear spatially invariant systems with sector nonlinearities and a road traffic network.
\end{abstract}


\section{Introduction}

Existing tools for controller synthesis do not scale to nowadays' complex large-scale systems. In large-scale vehicle platooning, for instance, classical distributed/decentralized control designs result in nonuniformity in the convergence rate of solutions; i.e., as the number of participating subsystems goes to infinity, the resulting network becomes unstable~\cite{Jovanovic.2005b}. Infinite networks, composed of infinitely many finite-dimensional subsystems, appear naturally as over-approximations of finite but very large networks with possibly unknown numbers of subsystems~\cite{Jovanovic.2005b}.%

Infinite networks appear in a wide variety of applications. Spatially invariant systems consisting of an infinite number of components interconnected to each other in the same pattern are studied in~\cite{BPD02,CIZ09} together with applications to, e.g., vehicle platoon formation~\cite{BeJ17}. Infinite systems also appear as representations of the solutions of linear and nonlinear partial differential equations over Hilbert spaces in terms of series expansions with respect to orthonormal or Riesz bases, see e.g.~\cite{LhS18}. A closely related approach relies on approximations of the system dynamics by partial differential and difference equations~\cite{Meurer.2013,kim2008pde}, which relies on a continuum approximation in space or in time and is particularly useful for consensus or coverage type problems. In addition, a closely related field is the ensemble control~\cite{Li11}, where the key objective is a simultaneous control of an infinite (and often uncountable) number of systems (neurons in the brain, flocks of birds, ensembles of quantum systems of the order of Avogadro's number $6 \cdot 10^{23}$) by a control signal which is the same for all subsystems.%

Most of the results on stability of infinite networks is devoted either to spatially invariant or to linear systems. Recently, several attempts have been made to relax such strong restrictions~\cite{DaP19,DMS19a,MKG20}, by introducing max-form small-gain theorems for infinite networks, where each subsystem is individually input-to-state stable (ISS)~\cite{Sontag.1989}. In~\cite{DaP19} it is shown that a countably infinite network of continuous-time ISS systems is ISS, provided that the gain functions capturing the influence from the neighboring subsystems are all less than identity which is conservative. By means of examples, it is shown in~\cite{DMS19a} that classic max-form small-gain conditions (SGCs) developed for finite-dimensional systems~\cite{Dashkovskiy.2010} do not guarantee the stability of infinite networks of ISS systems, even if all the systems are linear. To address this issue, more restrictive robust strong SGCs are developed in~\cite{DMS19a}.%

While the small-gain theorems in~\cite{DaP19,DMS19a} are formulated in terms of ISS Lyapunov functions, a trajectory-based nonlinear small-gain theorem for infinite networks is provided in~\cite{MKG20}. The main results in \cite{MKG20} show that if all subsystems are ISS with a uniform $\KL$-transient bound and with the gain operator satisfying the so-called monotone limit property, then the network is ISS. If the gain operator is linear, this condition is equivalent to the spectral small-gain condition (employed also in this work), and to the so-called uniform small-gain conditions.%

The starting point of the  ISS small-gain theory is to decompose a large-scale system (possibly an infinite network) into smaller subsystems and then analyze the stability properties of each subsystem individually.
In this way, it is assumed that each subsystem is ISS with respect to the neighboring subsystems, i.e., the inputs from other subsystems act as disturbances. Then, if the influence of the subsystems on each other is small enough, which is mathematically described by a SGC, stability of the overall system can be concluded. Small-gain theorems for finite-dimensional continuous-time systems can be found in~\cite{JMW96,DRW07,Dashkovskiy.2010}.%

For networks consisting of a countably infinite number of finite-dimensional continuous-time systems, we develop a SGC, which is tight in the sense that it cannot be relaxed under the assumptions we make on the network. We assume that each subsystem is exponentially ISS with respect to internal and external inputs and equipped with an exponential ISS Lyapunov function. The associated gain functions reflecting the interaction with neighbors are assumed to be linear. Such a scenario leads to several nontrivialities. In particular, the gain operator, which collects all the information about the internal gains, acts in an infinite-dimensional space, in contrast to couplings of just $N\in\N$ systems of arbitrary nature. This calls for a careful choice of an infinite-dimensional state space of the overall network, and 	motivates the use of the theory of positive operators on ordered Banach spaces for the small-gain analysis. 

\emph{We establish that if the gain operator $\Psi$, which is a positive operator, has spectral radius less than one:
\begin{equation}\label{eq:Spectral-SGC-Intro}
  r(\Psi) < 1,%
\end{equation}
}\emph{then the whole interconnection is exponentially ISS}. Furthermore, in our main result (cf.\ Theorem~\ref{MT}), we construct a so-called \emph{coercive} exponential ISS Lyapunov function for the overall network.%

Our main result is a \emph{nontrivial} generalization of Proposition 3.3 in~\cite{Dashkovskiy.2011b} from finite networks to infinite networks. The result in~\cite{Dashkovskiy.2011b} basically relies on \cite[Lem.~3.1]{Dashkovskiy.2011b}, which is a consequence of the Perron-Frobenius theorem. However, existing infinite-dimensional versions of the Perron-Frobenius theorem including the Krein-Rutman theorem \cite{KreiRut48}, are \emph{not} applicable to our setting as they require at least quasi-compactness of the gain operator, which is a quite strong assumption. We derive a technical lemma (cf.~Lemma~\ref{lem_smallgain}) which shows that under certain conditions there exists an infinite vector of scaling coefficients which is instrumental for the construction of a coercive exponential ISS Lyapunov function for the overall system. For classical results in the theory of positive operators in Banach spaces we refer the interested readers to \cite{K59}.%

The proposed small-gain criterion for the stability analysis of the network in \eqref{eq:Spectral-SGC-Intro} has several powerful characterizations in terms of the monotone bounded invertibility property and uniform small-gain conditions, as developed in \cite{MKG20, GlM20} and briefly explained in Section~\ref{sec:Verification of the spectral small-gain condition}, which help us to check this condition in practice.%

We illustrate the effectiveness of our results by applying them to \emph{nonlinear spatially invariant systems} with sector nonlinearities and a \emph{road traffic network}. We also add that the small-gain theorem shown in this paper has been extended recently in~\cite{NMK20b} to the case of ISS with respect to closed sets, and applied to study stability of infinite \emph{time-varying} networks, to consensus in \emph{infinite}-agent systems, as well as to the design of \emph{distributed observers} for infinite networks.%

The work in~\cite{DMS19a} is close in spirit to the present work, since in both the stability of the network is studied on the basis of the knowledge of ISS Lyapunov functions for the subsystems and the knowledge of the gain structure. However, in~\cite{DMS19a}, the ISS Lyapunov functions for the subsystems are defined in an implication form and the gain operator is used in a max formulation, which makes it \emph{nonlinear}, even if all the gains are linear. In contrast to~\cite{DMS19a}, in the present work we assume the existence of exponential ISS Lyapunov functions for the subsystems in a dissipative form and assume that the gain operator is defined in a sum form. These differences make the results of this paper and the methods employed in our analysis quite different from those of~\cite{DMS19a}. Furthermore, our tight small-gain criterion \eqref{eq:Spectral-SGC-Intro} for the stability analysis of the network is tailored to the exponential ISS Lyapunov functions which in this case is easier to check than the nonlinear small-gain conditions in~\cite{DMS19a}. In addition, the coercivity of the constructed Lyapunov function ensures a uniform decay rate of solutions of the network.%

This paper is organized as follows: First, relevant notation and the problem statement are given in Section~\ref{sec:Notation-and-problem-statement}. In Section~\ref{sec:Well-posedness}, we briefly discuss well-posedness of infinite networks. The notion of exponential ISS for infinite-dimensional systems in a Banach space and related Lyapunov properties are presented in Section~\ref{sec:ISS}. The technical results on the gain operator are made precise in Section~\ref{sec:The-gain-operator-and-its-properties}. The main small-gain theorem is presented in Section~\ref{sec:Small-gain-theorem}. The effectiveness of our results is verified through several examples in Section~\ref{sec:Examples}. Section~\ref{sec:Conclusions} concludes the paper.%

\section{Notation and Problem Statement}\label{sec:Notation-and-problem-statement}

\subsection{Notation}

We write $\N = \{1,2,3,\ldots\}$ for the set of positive integers. $\R$ denotes the reals and $\R_+ := \{t \in \R : t \geq 0\}$ the nonnegative reals. For vector norms on finite- and infinite-dimensional vector spaces, we write $|\cdot|$. For associated operator norms, we use the notation $\|\cdot\|$. We write $A\trn$ for the transpose of a matrix $A$ (which can be finite or infinite). We typically use Greek letters for infinite matrices and Latin ones for finite matrices. Elements of $\R^n$ are by default regarded as column vectors and we write $x\trn \cdot y$ for the Euclidean inner product of two vectors $x,y \in \R^n$. We use the same notation for dot products of vectors with infinitely many components. If $V:\R^n \rightarrow \R$ is a differentiable function, we write $\nabla V(x)$ for its gradient at $x$, which is a row vector by convention. By $\ell^p$, $p\in[1,\infty]$, we denote the Banach space of all real sequences $x = (x_i)_{i\in\N}$ with finite $\ell^p$-norm $|x|_p<\infty$, where $|x|_p = (\sum_{i=1}^{\infty}|x_i|^p)^{1/p}$ for $p < \infty$ and $|x|_{\infty} = \sup_{i\in\N}|x_i|$. We write $L^{\infty}(\R_+,\R^n)$ for the Banach space of essentially bounded measurable functions from $\R_+$ to $\R^n$. If $X$ is a Banach space, we write $r(T)$ for the spectral radius of a bounded linear operator $T:X \rightarrow X$ and $L(X)$ for the space of all bounded linear operators on $X$. The notation $C^0(X,Y)$ stands for the set of all continuous mappings $f:X \rightarrow Y$ between metric spaces $X$ and $Y$. Given a metric space $X$, we write $\inner\, A$ for the interior of a subset $A \subset X$. The right upper Dini derivative of a function $\gamma:\R \rightarrow \R$ at $t\in\R$ is defined by%
\begin{equation*}
  \rmD^+\gamma(t) := \limsup_{h \rightarrow 0+}\frac{1}{h}\big(\gamma(t + h) - \gamma(t)\big),%
\end{equation*}
and is allowed to assume the values $\pm \infty$. Analogously, the right lower Dini derivative of $\gamma$ at $t$ is defined by%
\begin{equation*}
  \rmD_+\gamma(t) := \liminf_{h \rightarrow 0+}\frac{1}{h}\big(\gamma(t + h) - \gamma(t)\big).%
\end{equation*}
Finally, we introduce the following classes of comparison functions which are frequently used in Lyapunov stability theory.%
\begin{align*}
  \PC &:= \bigl\{\gamma \in C^0(\R_+,\R_+) : \gamma(0)=0,\gamma(r)>0,\ \forall r > 0 \bigr\},\\
    \KC &:= \bigl\{\gamma \in \PC : \gamma \mbox{ is strictly increasing} \bigr\},\\
		\KC_{\infty} &:= \bigl\{\gamma \in \KC : \lim_{t\rightarrow\infty}\gamma(t)=\infty \bigr\},\\
		\LC &:= \bigl\{\gamma \in C^0(\R_+,\R_+) : \gamma \mbox{ is strictly decreasing with }\\
		&\qquad\qquad\qquad\qquad\qquad\qquad\qquad\qquad \lim_{t\rightarrow\infty}\gamma(t)=0 \bigr\},\\
  \KC\LC &:= \bigl\{\beta \in C^0(\R_+ \tm \R_+,\R_+) : \beta(\cdot,t) \in \KC\ \forall t\geq0,\\
	&\qquad\qquad\qquad\qquad\qquad\qquad\qquad \beta(r,\cdot) \in \LC\ \forall r>0 \bigr\}.%
\end{align*}

\subsection{Infinite interconnections}%

We study the interconnection of countably many systems, each given by a finite-dimensional ordinary differential equation (ODE). Using $\N$ as the index set (by default), the $i$-th subsystem is written as%
\begin{equation}\label{eq_ith_subsystem}
  \Sigma_i:\quad \dot{x}_i = f_i(x_i,\bar{x}_i,u_i).%
\end{equation}
The family $(\Sigma_i)_{i\in\N}$ comes together with sequences $(n_i)_{i\in\N}$, $(m_i)_{i\in\N}$ of positive integers and finite index sets $I_i \subset \N \backslash \{i\}$, $i\in\N$, so that the following assumptions hold.
\begin{itemize}
\item The state vector $x_i$ of $\Sigma_i$ is an element of $\R^{n_i}$.%
\item The vector $\bar{x}_i$ is composed of the state vectors $x_j$, $j \in I_i$. The order of these vectors plays no particular role (as the index set $\N$ does not), so we do not specify it.%
\item The external input vector $u_i$ is an element of $\R^{m_i}$.%
\item The right-hand side is a continuous function $f_i:\R^{n_i} \tm \R^{N_i} \tm \R^{m_i} \rightarrow \R^{n_i}$, where $N_i := \sum_{j\in I_i}n_j$.%
\item Unique local solutions of the ODE \eqref{eq_ith_subsystem} (in the sense of Carath\'eodory) exist for all initial states $x_{i0} \in \R^{n_i}$ and all locally essentially bounded functions $\bar{x}_i(\cdot)$ and $u_i(\cdot)$ (which are regarded as time-dependent inputs). We denote the corresponding solution by $\phi_i(t,x_{i0},\bar{x}_i,u_i)$.%
\end{itemize}
With the interconnection of the subsystems $\Sigma_i$ in mind, we interpret $\bar{x}_i(\cdot)$ as an \emph{internal input} and $u_i(\cdot)$ as an \emph{external input} (which may be a disturbance or a control input) for $\Sigma_i$ described in \eqref{eq_ith_subsystem}. The interpretation is that the subsystem $\Sigma_i$ is affected by finitely many neighbors, indexed by $I_i$, and its external input.%

To define the interconnection of the subsystems $\Sigma_i$, we consider the state vector $x = (x_i)_{i\in\N}$, the input vector $u = (u_i)_{i\in\N}$ and the right-hand side%
\begin{equation*}
  f(x,u) := (f_1(x_1,\bar{x}_1,u_1),f_2(x_2,\bar{x}_2,u_2),\ldots) \in \prod_{i\in\N}\R^{n_i}.%
\end{equation*}
The interconnection is then formally written as%
\begin{equation}\label{eq_interconnection}
  \Sigma:\quad \dot{x} = f(x,u).%
\end{equation}
To handle this infinite-dimensional ODE properly, we choose appropriate Banach spaces $X \subset \prod_{i\in\N}\R^{n_i}$ and $U \subset \prod_{i\in\N}\R^{m_i}$ and restrict $f$ to $X \tm U$.
As a natural choice, we use $\ell^p$-type spaces for both $X$ and $U$, and impose assumptions on $f$ to guarantee the existence and uniqueness of solutions. \emph{Our goal is then to show that $\Sigma$ is exponentially input-to-state stable (eISS) if each $\Sigma_i$ admits an eISS Lyapunov function and, additionally, a small-gain condition is satisfied.}%

\section{Well-Posedness}\label{sec:Well-posedness}

We want to model the state space $X$ of $\Sigma$ as a Banach space of sequences $x = (x_i)_{i\in\N}$ with $x_i \in \R^{n_i}$. The most natural choice is an $\ell^p$-space. To define such a space, we first fix a norm on each $\R^{n_i}$ (that might not only depend on the dimension $n_i$ but also on the index $i$). For brevity, we omit the index in our notation and simply write $|\cdot|$ for each of these norms. Then, for every $p \in [1,\infty)$, we put%
\begin{equation*}
  \ell^p(\N,(n_i)) := \Bigl\{ x = (x_i)_{i\in\N} : x_i \in \R^{n_i},\ \sum_{i\in\N}|x_i|^p < \infty \Bigr\}%
\end{equation*}
and equip this space with the norm $|x|_p := (\sum_{i\in\N}|x_i|^p)^{1/p}$.\footnote{The notation $\ell^p(\N,(n_i))$ should not be confused with $\ell^p$, which denotes the standard $\ell^p$-space of real sequences with finite $\ell^p$-norm.} Additionally, we introduce%
\begin{equation*}
  \ell^{\infty}(\N,(n_i)) := \Bigl\{ x = (x_i)_{i\in\N} : x_i \in \R^{n_i},\ \sup_{i\in\N}|x_i| < \infty \Bigr\},%
\end{equation*}
and equip this space with the norm $|x|_{\infty} := \sup_{i\in\N}|x_i|$. The following proposition is proved with standard arguments, see e.g.~\cite{DunSch57}. Hence, we omit the details.%

\begin{proposition}\label{prop1}
The following statements hold:%
\begin{enumerate}
\item[(a)] For each choice of norms on $\R^{n_i}$, $i\in\N$, and each $p \in [1,\infty]$, the associated space $\ell^p(\N,(n_i))$ equipped with the norm $|\cdot|_p$ is a Banach space.%
\item[(b)] For each $1 \leq p < \infty$, the Banach space $\ell^p(\N,(n_i))$ is separable.%
\item[(c)] For each pair $(p,q)$ with $1 \leq p < q \leq \infty$, the space $\ell^p(\N,(n_i))$ is continuously embedded into $\ell^q(\N,(n_i))$.%
\end{enumerate}
\end{proposition}

As the state space of the system $\Sigma$, we consider $X := \ell^p(\N,(n_i))$ for a fixed $p \in [1,\infty]$ (in the stability analysis, only finite $p$ will be considered). Similarly, for a fixed $q \in [1,\infty]$, we consider the \emph{external input space} $U := \ell^q(\N,(m_i))$, where we fix norms on $\R^{m_i}$ that we simply denote by $|\cdot|$ again. The space of admissible \emph{external input functions} is defined by\footnote{We use the letter $u$ both for elements of $U$ and $\UC$. Since it should become clear from the context if we refer to an input value or an input function, this should not lead to confusion.}%
\begin{align}\label{eq:Input-space}
  \UC := \bigl\{ u:\R_+ \rightarrow U : u & \mbox{ is strongly measurable}\nonumber\\
	 &\quad \mbox{and essentially bounded} \bigr\},%
\end{align}
where we recall that a strongly measurable function is defined as a (Borel) measurable function with a separable image. Since $\ell^q(\N,(m_i))$ is separable for all finite $q$, strong measurability reduces to measurability (and even to weak measurability, see \cite[Cor.~2, p.~73]{HiP00}) in all of these cases.%

A continuous mapping $\xi:J \rightarrow X$, defined on an interval $J = [0,T_*)$ with $T_* \in (0,\infty]$, is called a \emph{solution} of the infinite-dimensional ODE \eqref{eq_interconnection} with initial value $x^0 \in X$ for the external input $u\in\UC$ provided that $s \mapsto f(\xi(s),u(s))$ is an $X$-valued locally integrable function and%
\begin{equation*}
  \xi(t) = x^0 + \int_0^t f(\xi(s),u(s)) \rmd s%
\end{equation*}
holds for all $t \in J$, where the integral is the Bochner integral for Banach space valued functions. For a definition of the Bochner integral, we refer the interested readers to \cite{ABH11} and \cite[Ch.~III, Sec.~3.7]{HiP00}.%

If for each $x^0 \in X$ and $u \in \UC$ a unique local solution exists, we say that the system is \emph{well-posed} and write $\phi(\cdot,x^0,u)$ for any such solution. As usual, we consider the maximal extension of $\phi(\cdot,x^0,u)$ and write $J_{\max}(x^0,u)$ for its interval of existence. We say that the system is \emph{forward complete} if $J_{\max}(x^0,u) = \R_+$ for all $(x^0,u) \in X \tm \UC$.%

Denoting by $\pi_i:X \rightarrow \R^{n_i}$ the canonical projection onto the $i$-th component (this projection is a bounded linear operator) and writing $u(t) = (u_1(t),u_2(t),\ldots)$, we obtain%
\begin{align*}
 & \pi_i\phi(t,x^0,u) = x^0_i + \int_0^t \pi_i f\big(\phi(s,x^0,u),u(s)\big) \rmd s \\
	&= x^0_i + \int_0^t f_i(\pi_i\phi(s,x^0,u),(\pi_j\phi(s,x^0,u))_{j\in I_i},u_i(s)) \rmd s,%
\end{align*}
which implies that $t \mapsto \pi_i\phi(t,x^0,u)$ solves the ODE%
\begin{equation*}
  \dot{x}_i = f_i(x_i,\bar{x}_i,u_i)%
\end{equation*}
for the internal input $\bar{x}_i(\cdot) := (\pi_j\phi(\cdot,x^0,u))_{j\in I_i}$ and the external input $u_i(\cdot)$. Hence,%
\begin{equation*}
  \pi_i \phi(t,x^0,u) = \phi_i(t,x^0_i,\bar{x}_i,u_i) \mbox{\quad for all\ } t \in J_{\max}(x^0,u).%
\end{equation*}

Sufficient conditions for the existence and uniqueness of solutions (and forward completeness) can be obtained from the general theory of Carath\'eodory differential equations on Banach spaces, see \cite{AulWan96} as a general reference for systems with bounded generators.%

The following theorem provides a set of conditions which is sufficient for well-posedness and forward completeness.%

\begin{theorem}\label{thm_awa}
Assume that the system $\Sigma$ with state space $X = \ell^p(\N,(n_i))$ and external input space $U = \ell^q(\N,(m_i))$ for arbitrary $p,q \in [1,\infty]$ satisfies the following assumptions.
\begin{enumerate}
\item[(i)] $f(x,u) \in X$ for all $(x,u) \in X \tm U$.%
\item[(ii)] For every $u \in U$, the mapping $f(\cdot,u):X \rightarrow X$ is continuous.%
\item[(iii)] For every $x \in X$, the mapping $f(x,\cdot):U \rightarrow X$ is continuous.%
\item[(iv)] For each $u \in \UC$, there exist locally integrable functions $\ell,\ell_0:\R_+ \rightarrow \R_+$ such that%
\begin{subequations}
\label{eq:Well-posedness-iv}
\begin{align}
  |f(x^1,u(t)) - f(x^2,u(t))|_p &\leq \ell(t) |x^1 - x^2|_p, \label{eq:Well-posedness-iv-1} \\
	|f(0,u(t))|_p &\leq \ell_0(t) \label{eq:Well-posedness-iv-2}%
\end{align}
\end{subequations}
hold for almost all $t\in \R_+$ and all $x^1,x^2 \in X$.%
\end{enumerate}
Then for every initial value $x^0 \in X$ and every external input $u \in \UC$ a unique solution $\phi(\cdot,x^0,u):\R_+ \rightarrow X$ exists and for any $u\in\UC$, the mapping $(t,x^0) \mapsto \phi(t,x^0,u)$ is continuous on $\R_+ \tm X$.%
\end{theorem}

Observe that Assumption (iii) in Theorem~\ref{thm_awa} implies that the function $t \mapsto f(x,u(t)),\quad \R_+ \rightarrow X$ is strongly measurable for each $x \in X$ and $u \in \UC$, since it is the composition of the strongly measurable function $u$ and the continuous function $u \mapsto f(x,u)$. The theorem then follows immediately from \cite[Thm.~2.4]{AulWan96}.%

From this result on global existence and uniqueness of solutions we can easily deduce a result on local existence and uniqueness which requires less strict assumptions.%

\begin{corollary}
Assume that system $\Sigma$ with state space $X = \ell^p(\N,(n_i))$ and external input space $U = \ell^q(\N,(m_i))$ for arbitrary $p,q \in [1,\infty]$ satisfies assumptions (i)--(iii) of Theorem~\ref{thm_awa} as well as the following property:
\begin{enumerate}
\item[(iv)] For each $u \in \UC$ and $x^0 \in X$, there exist $\delta>0$ and locally integrable functions $\ell,\ell_0:\R_+ \rightarrow \R_+$ such that%
\begin{subequations}
\label{eq:Well-posedness-iv-local}
\begin{align}
  |f(x^1,u(t)) - f(x^2,u(t))|_p &\leq \ell(t) |x^1 - x^2|_p, \label{eq:Well-posedness-iv-1x} \\
	|f(x^0,u(t))|_p &\leq \ell_0(t), \label{eq:Well-posedness-iv-2x}%
\end{align}
\end{subequations}
hold for almost all $t\in \R_+$ and all $x^1,x^2 \in B_{\delta}(x^0)$.%
\end{enumerate}
Then for every initial value $x^0 \in X$ and every external input $u \in \UC$ a unique solution $\phi(\cdot,x^0,u):[0,t_{\max}(x^0,u)) \rightarrow X$ exists, where $0 < t_{\max}(x^0,u) \leq \infty$.%
\end{corollary}

\begin{IEEEproof}
Having fixed some $x^0 \in X$ and $u\in\UC$, let $\delta>0$ and $\ell,\ell_0$ be as in (iv) in Theorem \ref{thm_awa}. We modify the right-hand side $f$ so that it satisfies the assumptions of Theorem \ref{thm_awa}. To this end, we use the retraction $r_0:X \rightarrow X$ defined by%
\begin{equation*}
  r_0(x) := \left\{ \begin{array}{cc}
	                      x & \mbox{if } |x|_p \leq 1,\\
											  \frac{x}{|x|_p} & \mbox{if } |x|_p > 1.%
									\end{array} \right.%
\end{equation*}
It is not hard to show that $r_0$ satisfies a global Lipschitz condition, see, e.g., the proof of \cite[Thm.~5.1]{Sle89} for a similar argument. Now we define%
\begin{equation*}
  r_{x^0,\delta}(x) := \delta r_0\Bigl(\frac{x - x^0}{\delta}\Bigr) + x^0,%
\end{equation*}
and observe that $|x - x^0|_p \leq \delta$ implies $r_{x^0,\delta}(x) = x$, while $|x - x^0| > \delta$ implies $r_{x^0,\delta}(x) = \delta\frac{x - x_0}{|x - x^0|_p} + x^0 \in \overline{B_{\delta}(x^0)}$. Moreover, $r_{x^0,\delta}$ is globally Lipschitz continuous with Lipschitz constant $2$. We define $\tilde{f}(x,u) := f(r(x),u)$ and consider system%
\begin{equation*}
  \tilde{\Sigma}:\quad \dot{x} = \tilde{f}(x,u).%
\end{equation*}
Using the Lipschitz continuity of $r$, it is easy to verify that $\tilde{\Sigma}$ satisfies the assumptions of Theorem \ref{thm_awa}. Obviously, every solution of $\tilde{\Sigma}$ starting in $x^0$ is locally a solution of $\Sigma$. This proves the corollary.
\end{IEEEproof}

\begin{remark}\label{rem:InfODEs-Special-case-of-gen-systems} 
We note that well-posed systems \eqref{eq_interconnection} are control systems in the sense of \cite[Def.~1]{MiW18b}, and thus a number of results in the general ISS theory of infinite-dimensional systems \cite{MiP19} are valid for the system \eqref{eq_interconnection}.\hfill$\diamond$
\end{remark}

\begin{example}\label{ex_linear_general}
Assume that the subsystems $\Sigma_i$ are linear:%
\begin{equation*}
  \Sigma_i:\quad \dot{x}_i = A_ix_i + \tilde{A}_i \bar{x}_i + B_iu_i,%
\end{equation*}
with matrices $A_i \in \R^{n_i \tm n_i}$, $\tilde{A}_i \in \R^{n_i \tm N_i}$, and $B_i \in \R^{n_i \tm m_i}$. Given a choice of norms on $\R^{n_i}$, on $\R^{N_i} = \prod_{j\in I_i}\R^{n_j}$ consider the product norm%
\begin{equation*}
  |\bar{x}_i| := \sum_{j\in I_i}|x_j|.%
\end{equation*}

We make the following assumptions.%
\begin{itemize}
\item The operator norms (with respect to the chosen vector norms on $\R^{n_i}$, $\R^{N_i}$ and $\R^{m_i}$) of the linear operators $A_i$, $\tilde{A}_i$, and $B_i$ are uniformly bounded over $i\in\N$.%
\item $1 \leq p = q < \infty$.%
\item There exists $m\in\N$ such that $I_i \subset [i-m,i+m] \cap \N$ for all $i\in\N$.%
\end{itemize}
Now, we show that Assumptions (i)-(iv) in Theorem \ref{thm_awa} are satisfied. Take $x \in X = \ell^p(\N,(n_i))$ and $u \in U = \ell^p(\N,(m_i))$.%

By using the inequality $(a+b)^p \leq 2^{p-1}(a^p + b^p)$ (which holds for all $a,b\geq 0$ due to the convexity of $a \mapsto a^p$) repeatedly, we obtain%
\begin{align*}
  \sum_{i=1}^{\infty} |A_ix_i& + \tilde{A}_i\bar{x}_i + B_iu_i|^p \\
	&\leq C_1 \sum_{i=1}^{\infty}|x_i|^p + C_2\sum_{i=1}^{\infty}\sum_{j \in I_i}|x_j|^p + C_3\sum_{i=1}^{\infty}|u_i|^p\\
	&\leq C_1 |x|_p^p + C_2 (2m+1) |x|_p^p + C_3|u|_p^p < \infty,%
\end{align*}
where $C_1,C_2,C_3>0$ are appropriately chosen constants. This shows that Assumption (i) in Theorem \ref{thm_awa} is satisfied.%

To see that Assumption (ii) holds, observe that for any $x^1,x^2 \in X$ and $u \in U$ we have%
\begin{align}
  |f(x^1,u) -& f(x^2,u)|_p^p = \sum_{i=1}^{\infty} \bigl|A_i(x_i^1 - x_i^2) + \tilde{A}_i(\bar{x}_i^1 - \bar{x}_i^2)\bigr|^p\nonumber\\
	&\leq C_1 \sum_{i=1}^{\infty} |x_i^1 - x_i^2|^p + C_2 \sum_{i=1}^{\infty}|\bar{x}_i^1 - \bar{x}_i^2|^p\nonumber\\
	&\leq C  |x^1 - x^2|_p^p, \label{eq:Ass-iii-linear}%
\end{align}
for some constant $C>0$. Furthermore,%
\begin{equation*}
  |f(x,u^1) - f(x,u^2)|_p^p \leq C_3 \sum_{i=1}^{\infty}|u_i^1 - u_i^2|^p = C_3 |u^1 - u^2|_p^p,%
\end{equation*}
which implies that Assumption (iii) is satisfied.%

Finally, for Assumption (iv), note that the inequality \eqref{eq:Well-posedness-iv-1} is valid with a constant function $\ell$, due to \eqref{eq:Ass-iii-linear}. For \eqref{eq:Well-posedness-iv-2}, note that%
\begin{equation*}
  |f(0,u(t))|_p^p \leq C_3 \sum_{i=1}^{\infty} |u_i(t)|^p = C_3 |u(t)|_p^p.%
\end{equation*}  
By assumption, $u$ is essentially bounded as a function from $\R_+$ into $\ell^p(\N,(m_i))$, which implies that $|u(\cdot)|_p$ is locally integrable.\hfill$\square$
\end{example}

\section{Exponential Input-to-State Stability}\label{sec:ISS}

Having a well-posed interconnection \eqref{eq_interconnection} with state space $X = \ell^p(\N,(n_i))$ and external input space $U = \ell^q(\N,(m_i))$, it is natural to study its stability with respect to both initial conditions and external inputs. The concept of input-to-state stability is suitable for both of these purposes.%

We equip the (linear) space $\UC$ of external inputs with the sup-norm%
\begin{equation*}
  |u|_{q,\infty} := \esssup_{t\geq0}|u(t)|_q%
\end{equation*}
and work with the following definition of input-to-state stability (cf.~\cite[Def.~6]{DaM13}).%

\begin{definition}
The system $\Sigma$ is called \emph{input-to-state stable (ISS)}  if it is forward complete and there exist functions $\beta \in \KC\LC$ and $\gamma \in \KC$ such that for any initial state $x^0 \in X$ and any $u\in\UC$ the corresponding solution satisfies%
\begin{equation*}
  |\phi(t,x^0,u)|_p \leq \beta(|x^0|_p,t) + \gamma(|u|_{q,\infty}) \mbox{\quad for all\ } t \geq 0.%
\end{equation*}
\end{definition}

If the decay of the norm of $\phi(t,x^0,u)$ is exponential in $t$, the system is called \emph{exponentially input-to-state stable}. The precise definition reads as follows.%

\begin{definition}\label{def:eISS}
The system $\Sigma$ is called \emph{exponentially input-to-state stable (eISS)} if it is forward complete and there are constants $a,M>0$ and $\gamma \in \KC$ such that for any initial state $x^0 \in X$ and any $u\in\UC$ the corresponding solution satisfies%
\begin{equation*}
  |\phi(t,x^0,u)|_p \leq M\rme^{-at}|x^0|_p + \gamma(|u|_{q,\infty}) \mbox{\quad for all\ } t \geq 0.%
\end{equation*}
\end{definition}

We note that, by the causality of $\Sigma$, eISS implies the following inequality:%
\begin{equation*}
  |\phi(t,x^0,u)|_p \leq M\rme^{-at}|x^0|_p + \gamma\bigl(\esssup_{0 \leq s \leq t}|u(s)|_q\bigr).%
\end{equation*}

For any continuous function $V:X \rightarrow \R$, let us define the \emph{orbital derivative} at $x \in X$ for the external input $u \in \UC$ by%
\begin{equation*}
  \rmD^+ V_u(x) := \rmD^+ V(\phi(t,x,u))_{|t=0},%
\end{equation*}
where the right-hand side is the right upper Dini derivative of the function $t \mapsto V(\phi(t,x,u))$, evaluated at $t=0$.%

Exponential input-to-state stability is implied by the existence of an exponential ISS Lyapunov function, which we define in a dissipative form as follows.%

\begin{definition}\label{def:eISS-Lyapunov-Function}
A continuous function $V:X \rightarrow \R_+$ is called an \emph{eISS Lyapunov function} for the system $\Sigma$ if there exist constants $\underline{\omega},\overline{\omega},b,\kappa>0$ and $\gamma \in \KC_{\infty}$ such that%
\begin{subequations}\label{eq_eiss_lyap_props}
\begin{align}
  \underline{\omega} |x|_p^b &\leq V(x) \leq \overline{\omega} |x|_p^b, \label{eq_eiss_lyap_props-1}\\
  \rmD^+ V_u(x) &\leq -\kappa V(x) + \gamma(|u|_{q,\infty}) \label{eq_eiss_lyap_props-2}%
\end{align}
\end{subequations}
hold for all $x \in X$ and $u \in \UC$.%
\end{definition}

ISS Lyapunov functions, for which there is $\psi\in\KC_\infty$ such that $V(x) \geq \psi(|x|)$ for all $x \in X$, are frequently called \emph{coercive} ISS Lyapunov functions, as opposed to \emph{non-coercive} ISS Lyapunov functions, introduced in \cite{MiW18b}, for which a weaker estimate $V(x)>0$ for all $x\neq 0$, holds.%

The importance of eISS Lyapunov functions is due to the following result (the proof is a variation of the arguments given in \cite[Thm.~1]{DaM13}, and thus only some of the steps are provided).%

\begin{proposition}\label{prop_eiss}
If there exists an eISS Lyapunov function for $\Sigma$, then $\Sigma$ is eISS.%
\end{proposition}

\begin{IEEEproof}
Let $V$ be an eISS Lyapunov function as in Definition~\ref{def:eISS-Lyapunov-Function} with corresponding constants $\underline{\omega},\overline{\omega},b,\kappa>0$ and a function $\gamma \in \KC_{\infty}$.%

Pick any $\varepsilon \in (0,\kappa)$ and define $\chi(r) := \frac{1}{\kappa-\varepsilon} \gamma(r)$. For all $x \in X$ and $u\in\UC$, we obtain%
\begin{equation*}
  V(x) \geq \chi(|u|_{q,\infty}) \quad\Rightarrow\quad \rmD^+ V_u(x) \leq -\varepsilon V(x).%
\end{equation*}
By arguing similarly to \cite[Thm.~1]{DaM13}, we obtain the following inequality for all $x \in X$, $u\in\UC$, and $t\geq 0$:%
\begin{equation*}
  V(\phi(t,x,u)) \leq \rme^{-\varepsilon t} V(x) + \chi(|u|_{q,\infty}).%
\end{equation*}
In view of \eqref{eq_eiss_lyap_props-1}, we get%
\begin{equation*}
  \underline{\omega} |\phi(t,x,u)|_p^b \leq \rme^{-\varepsilon t} \overline{\omega} |x|_p^b + \chi(|u|_{q,\infty}),%
\end{equation*}
which, by the monotonicity of $\gamma$, implies that%
\begin{align*}
 |\phi(t,x,u)|_p 
&\leq \Big(\rme^{-\varepsilon t} \frac{\overline{\omega}}{\underline{\omega}} |x|_p^b + \chi(|u|_{q,\infty})\Big)^{\frac{1}{b}}\\
&\leq \Big(2 \frac{\overline{\omega}}{\underline{\omega}}\Big)^{\frac{1}{b}} \rme^{-\frac{\varepsilon}{b} t} |x|_p + \Big(2\chi(|u|_{q,\infty})\Big)^{\frac{1}{b}},%
\end{align*}
showing the exponential ISS property for $\Sigma$.%
\end{IEEEproof}

\begin{remark}\label{rem:variations-of-eISS} 
One can find some variations of the eISS property in \cite[p.~2736]{HLT08}. Observe that in our definition of eISS Lyapunov functions, in addition to the exponential decay along the trajectory, we also assume that bounds of the form \eqref{eq_eiss_lyap_props-1} hold. The inequalities \eqref{eq_eiss_lyap_props-1} are needed to ensure that the existence of an eISS Lyapunov function implies eISS. The exponential decay of a Lyapunov function along trajectories alone is not sufficient for this implication.\hfill$\diamond$
\end{remark}

\section{The Gain Operator and its Properties}\label{sec:The-gain-operator-and-its-properties}

Our main objective is to develop conditions for input-to-state stability of the interconnection of countably many subsystems \eqref{eq_ith_subsystem}, depending on certain stability properties of the subsystems.%

\subsection{Assumptions on the subsystems}

We assume that each subsystem $\Sigma_i$, given by \eqref{eq_ith_subsystem}, is  exponentially ISS and there exist continuous eISS Lyapunov functions with linear gains for all $\Sigma_i$. Restating the concept of an eISS Lyapunov function (Definition~\ref{def:eISS-Lyapunov-Function}) for the subsystem $\Sigma_i$, we see that the gain $\gamma$ in this definition indicates the influence of the aggregated input onto the system. For our purposes, this information is not sufficient as we would like to know how each $j$-th subsystem influences each $i$-th subsystem as in the next assumption.%

\begin{assumption}\label{ass_vi_existence}
For each $i\in\N$ there exists a continuous function $V_i:\R^{n_i} \rightarrow \R_+$, satisfying for certain $p,q\in[1,\infty)$ the following properties.%
\begin{itemize}
\item[(i)] There are constants $\underline{\alpha}_i,\overline{\alpha}_i>0$ so that for all $x_i \in \R^{n_i}$%
\begin{equation}\label{eq_viest}
  \underline{\alpha}_i|x_i|^p \leq V_i(x_i) \leq \overline{\alpha}_i|x_i|^p.%
\end{equation}
\item[(ii)] There are constants $\lambda_i,\gamma_{ij},\gamma_{iu}>0$ so that the following holds: for each $x_i\in\R^{n_i}$,  $u_i \in L^{\infty}(\R_+,\R^{m_i})$ and each internal input $\bar{x}_i \in C^0(\R_+,\R^{N_i})$ and for almost all $t$ in the maximal interval of existence of $\phi_i(t) := \phi_i(t,x_i,\bar{x}_i,u_i)$ one has%
\begin{align}\label{eq_nablaviest}									
\begin{split}
\hspace{-4mm}\rmD^+ (V_i \circ \phi_i)(t) \leq -&\lambda_i V_i(\phi_i(t)) + \sum_{j \in I_i} \gamma_{ij}V_j(x_j(t))\\
&\qquad\qquad\qquad + \gamma_{iu}|u_i(t)|^q,%
\end{split}
\end{align}
where we denote the components of $\bar{x}_i$ by $x_j(\cdot)$.%
\item[(iii)] For all $t$ in the maximal interval of the existence of $\phi_i$%{\rm :}%
\begin{equation*}
  \rmD_+(V_i \circ \phi_i)(t) < \infty.%
\end{equation*}
\end{itemize}
\end{assumption}

\begin{remark}\label{rem:Dissipation-conditions-C1-Vi} 
Often eISS Lyapunov functions $V_i$ for subsystems \eqref{eq_ith_subsystem} are assumed to be continuously differentiable. In this case, the dissipative conditions \eqref{eq_nablaviest} for ISS Lyapunov functions $V_i$ can be formulated in a computationally simpler style:%
\begin{align}\label{eq_nablaviest-C1-case}
\begin{split}
  &\nabla V_i(x_i) \cdot f_i(x_i,\bar{x}_i,u_i)\\
	&\qquad\qquad\leq -\lambda_i V_i(x_i) + \sum_{j \in I_i}\gamma_{ij}V_j(x_j) + \gamma_{iu}|u_i|^q.%
\end{split}
\end{align}
These conditions have to be valid for all $x_i\in\R^{n_i}$, $u_i\in \R^{m_i}$ and $\bar{x}_i \in \R^{N_i}$. The expression on the left-hand side of inequality \eqref{eq_nablaviest-C1-case} represents a formula for the computation of the orbital derivative of $V_i$ under the assumption that $V_i$ is smooth enough. In this case, item (iii) in Assumption \ref{ass_vi_existence} is automatically fulfilled. The proof of the corresponding small-gain theorem goes along the same lines as in Theorem~\ref{MT}. \hfill$\diamond$%
\end{remark}

We furthermore assume that the following uniformity conditions hold for the constants introduced above, which turn out to be crucial for the construction of a coercive eISS Lyapunov function for the interconnection of the subsystems $\Sigma_i$.%

\begin{assumption}\label{ass_external_gains}
\begin{enumerate}
\item[(a)] There are constants $\underline{\alpha},\overline{\alpha}>0$ so that for all $i\in \N$%{\rm :}%
\begin{equation}\label{eq_uniformity_alpha}
  \underline{\alpha} \leq \underline{\alpha}_i \leq \overline{\alpha}_i \leq \overline{\alpha}.%
\end{equation}
\item[(b)] There is a constant $\underline{\lambda}>0$ so that for all $i\in\N$
\begin{equation}\label{eq_uniformity_lambda}
  \underline{\lambda} \leq \lambda_i.%
\end{equation}
\item[(c)] There is a constant $\overline{\gamma}_u>0$ so that for all $i\in\N$%
\begin{equation}\label{eq_uniformity_gammaiu}
  \gamma_{iu} \leq \overline{\gamma}_u.%
\end{equation}
\end{enumerate}
\end{assumption}

Assumption \ref{ass_vi_existence} enforces stability properties of the subsystems $\Sigma_i$. In order to speak about the interconnection of all subsystems in \eqref{eq_ith_subsystem}, we should define the state space for the interconnection as well as the space of input values. The inequalities \eqref{eq_viest} and \eqref{eq_nablaviest} suggest the following choice: 
\begin{center}
$X = \ell^p(\N,(n_i))$ and $U = \ell^q(\N,(m_i))$.%
\end{center}

We thus make the following well-posedness assumption.%

\begin{assumption}\label{ass_wellposedness}
The system $\Sigma$ with state space $X = \ell^p(\N,(n_i))$ and external input space $U = \ell^q(\N,(m_i))$ is well-posed.%
\end{assumption}

\begin{remark}\label{rem:Choice-of-a-state space}
We define the state and input space for the overall interconnected system based on the values of the parameters $p,q$ which we obtain from Assumption \ref{ass_vi_existence}. However, the choice of the state space depends also on the physical meaning of the variables $x_i$. For example, if $x_i$ represents a mass, and we are interested in the dynamics of the total mass of a system, then a reasonable choice of the state space for the interconnection is $X = \ell^1(\N,(n_i))$, and if $x_i$ represents a velocity and we are interested in the dynamics of the total kinetic energy of the system, then it is natural to choose $X = \ell^2(\N,(n_i))$. Therefore, to meet the needs of applications, one should construct the ISS Lyapunov functions $V_i$ for some specific values of $p,q$. \hfill$\diamond$%
\end{remark}

We note that inequalities \eqref{eq_viest} in terms of the $\KC_{\infty}$-functions $r \mapsto \underline{\alpha}r^p$ and $r \mapsto \overline{\alpha}r^p$ turn out to be crucial for a sum-type construction of an eISS Lyapunov function for $\Sigma$.%

In order to formulate a small-gain condition, we further introduce infinite nonnegative matrices by collecting the coefficients from \eqref{eq_nablaviest} as follows:%
\begin{equation*}
  \Lambda := \diag(\lambda_1,\lambda_2,\lambda_3,\ldots), \quad \Gamma := (\gamma_{ij})_{i,j\in\N},%
\end{equation*}
where we put $\gamma_{ij} := 0$ whenever $j \notin I_i$. We also introduce the infinite matrix%
\begin{equation}\label{eq:operator-A}
  \Psi := \Lambda^{-1}\Gamma = (\psi_{ij})_{i,j\in\N},\quad \psi_{ij} = \frac{\gamma_{ij}}{\lambda_i}.
\end{equation}

Under an appropriate boundedness assumption, the matrix $\Psi$ acts as a linear operator on $\ell^1$ by%
\begin{equation*}
  (\Psi x)_i = \sum_{j=1}^{\infty} \psi_{ij} x_j \mbox{\quad for all\ } i \in \N.%
\end{equation*}

We call $\Psi:\ell^1 \rightarrow \ell^1$ the \emph{gain operator} associated with the decay rates $\lambda_i$ and coefficients $\gamma_{ij}$.%

We make the following assumption, which is equivalent to $\Gamma$ being a bounded operator from $\ell^1$ to $\ell^1$.%

\begin{assumption}\label{ass_A_bounded}
The matrix $\Gamma = (\gamma_{ij})$ satisfies%
\begin{equation}\label{eq_Gamma_bounded}
  \|\Gamma\|_{1,1} = \sup_{j \in \N} \sum_{i=1}^{\infty} \gamma_{ij} < \infty,%
\end{equation}
where the double index on the left-hand side indicates that we consider the operator norm induced by the $\ell^1$-norm both on the domain and codomain of the operator $\Gamma$.
\end{assumption}

\begin{remark}\label{rem:uniformity_gamma} 
Assumption~\ref{ass_A_bounded} implies that there is a constant $\overline{\gamma}>0$ such that $ 0 < \gamma_{ij} \leq \overline{\gamma}$ for all $i\in\N$ and $j \in I_i$. \hfill$\diamond$%
\end{remark}

Under Assumptions~\ref{ass_A_bounded} and~\ref{ass_external_gains}(b), the gain operator $\Psi$ is bounded.%

\begin{lemma}\label{lem:Gain-Operator-Bounded} 
Suppose that Assumptions~\ref{ass_A_bounded} and~\ref{ass_external_gains}(b) hold. Then $\Psi:\ell^1 \rightarrow \ell^1$, defined by~\eqref{eq:operator-A}, is a bounded operator.%
\end{lemma}

\begin{IEEEproof}
It holds that 
\begin{equation*}
  \|\Psi\|_{1,1} = \sup_{j \in \N} \sum_{i=1}^{\infty} \psi_{ij}
	= \sup_{j \in \N} \sum_{i=1}^{\infty} \frac{\gamma_{ij}}{\lambda_{i}} 
	\leq \frac{1}{\underline{\lambda}}  \sup_{j \in \N} \sum_{i=1}^{\infty} \gamma_{ij} < \infty,%
\end{equation*}
which is equivalent to boundedness of $\Psi$.
\end{IEEEproof}

A sufficient (though not necessary) condition for \eqref{eq_Gamma_bounded} is provided by the following lemma. The proof is simple and is omitted here.%

\begin{lemma}\label{lem_Gamma_bounded}
If there exists $m\in\N$ so that $I_i \subset [i-m,i+m] \cap \N$ for all $i\in\N$ and $\gamma_{ij} \leq \overline{\gamma}$ for all $i,j\in\N$ with a constant $\overline{\gamma}>0$, then \eqref{eq_Gamma_bounded} holds.%
\end{lemma}

\subsection{Spectral radius of the gain operator}%

In this subsection, we prove an auxiliary result which yields the existence of an infinite vector $\mu \in \ell^{\infty}$ that can be used to construct an eISS Lyapunov function for $\Sigma$ from the individual Lyapunov functions $V_i$, under the assumption that $r(\Psi) < 1$ (the small-gain condition) holds for the spectral radius of $\Psi$.%

For an overview of the concepts from functional analysis used here, see Appendix \ref{sec_positive_operators}.%

The adjoint operator of $\Psi$ acts on $\ell^{\infty}$ (which is canonically identified with the dual space $(\ell^1)^*$) and can be described by the transpose $\Theta:=\Psi\trn$, $\Theta = (\theta_{ij}) = (\psi_{ji})$ as%
\begin{equation*}
  (\Theta x)_i = \sum_{j=1}^{\infty} \theta_{ij}x_j = \sum_{j=1}^{\infty} \frac{\gamma_{ji}}{\lambda_j}x_j \quad\forall x \in \ell^{\infty}.%
\end{equation*}

On the Banach space $\ell^{\infty}$, we consider the standard positive cone%
\begin{equation*}
  K := \left\{ (x_i)_{i\in\N} \in \ell^{\infty} : x_i \geq 0,\ \forall i \in \N \right\},%
\end{equation*}
and observe that the interior of $K$ is nonempty and given by%
\begin{equation*}
  \inner\, K = \left\{ x \in \ell^{\infty} : \exists \underline{x}>0 \mbox{ s.t. } x_i \geq \underline{x},\ \forall i\in\N \right\}.%
\end{equation*}

Clearly, $\Theta$ maps the cone $K$ into itself, hence, is a positive operator with respect to this cone. The partial order on $\ell^{\infty}$, induced by $K$, is given by%
\begin{equation*}
  x \leq y \quad \Leftrightarrow \quad x_i \leq y_i,\ \forall i \in \N.%
\end{equation*}

The next technical result will be used to prove the main result of the paper (cf.~Theorem \ref{MT}).%

\begin{lemma}\label{lem_smallgain}
Assume that the constants $\lambda_i$ are uniformly bounded from above. Then $r(\Theta)<1$ implies the following statements:%
\begin{itemize}
\item[(i)] there exist a vector $\mu = (\mu_i)_{i\in\N} \in \inner\, K$ and a constant $\lambda_{\infty}>0$ so that%
\begin{equation*}
  \frac{[\mu\trn(-\Lambda + \Gamma)]_i}{\mu_i} \leq -\lambda_{\infty} \mbox{\quad for all\ } i \in \N;%
\end{equation*}
\item[(ii)] for every $\rho>0$, we can choose the vector $\mu$ and the constant $\lambda_{\infty}$ in statement (i) in such a way that the following inequality holds:%
\begin{equation*}
  \lambda_{\infty} \geq (1 - r(\Theta))\underline{\lambda} - \rho.%
\end{equation*}
\end{itemize}
\end{lemma}

\begin{IEEEproof}
From Lemma~\ref{lem:rT_bound-new} and the assumption that $r(\Theta)<1$, it follows that there exist $\eta \in \inner\, K$ and $\tilde{r} \in (0,1)$ with $\Theta \eta \leq \tilde{r}\eta$, or equivalently $\eta\trn \Lambda^{-1}\Gamma \le \tilde{r} \eta\trn$. Defining $\mu\trn := \eta\trn\Lambda^{-1}$, the previous inequality can be rewritten as $\mu\trn (-\tilde{r}\Lambda + \Gamma) \le 0$, which we further transform into%
\begin{equation*}
  \mu\trn(-\Lambda + \Gamma) \le \mu\trn(-(1 - \tilde{r})\Lambda).%
\end{equation*}
Thus, we obtain%
\begin{equation*}
  \left[\mu\trn(-\Lambda + \Gamma)\right]_i \leq -(1 - \tilde{r}) \lambda_i \mu_i \leq -(1 - \tilde{r}) \underline{\lambda}\,\mu_i% 
\end{equation*}
with $\underline{\lambda}$ as in \eqref{eq_uniformity_lambda}. To show that $\mu \in \inner\, K$, we use the assumption that there exists an upper bound $\overline{\lambda}$ for the constants $\lambda_i$, which guarantees that the components of $\mu$ satisfy $\mu_i = \lambda_i^{-1} \eta_i \geq \overline{\lambda}^{-1}\eta_i$. Thus, statement (i) is proved.%

Statement (ii) follows from the fact that we can choose $\tilde{r}$ arbitrarily close to $r(\Theta)$.%
\end{IEEEproof}

\section{Small-Gain Theorem}\label{sec:Small-gain-theorem}

In this section, we prove that the interconnected system $\Sigma$ is exponentially ISS under the given assumptions, provided that the spectral radius of the gain operator satisfies $r(\Psi)<1$.

By Proposition \ref{prop_eiss}, our objective is reduced to finding an eISS Lyapunov function for the interconnection $\Sigma$ under the small-gain condition $r(\Psi)<1$. This is accomplished by the following \emph{small-gain theorem}, which is the main result of the paper.%

\begin{theorem}\label{MT}
Consider the infinite interconnection $\Sigma$, composed of subsystems $\Sigma_i$, $i\in\N$, with fixed $p,q \in [1,\infty)$, and let the following assumptions be satisfied.%
\begin{enumerate}
\item[(i)] $\Sigma$ is well-posed as a system with state space $X = \ell^p(\N,(n_i))$, space of input values $U = \ell^q(\N,(m_i))$, and the external input space $\UC$, as defined in \eqref{eq:Input-space}.%
\item[(ii)] Each $\Sigma_i$ admits a continuous eISS Lyapunov function $V_i$ so that Assumptions \ref{ass_vi_existence} and \ref{ass_external_gains} are satisfied.%
\item[(iii)] The operator $\Gamma:\ell^1 \rightarrow \ell^1$ is bounded, i.e., Assumption \ref{ass_A_bounded} holds.%
\item[(iv)] The spectral radius of $\Psi$ satisfies $r(\Psi) < 1$.%
\end{enumerate}
Then $\Sigma$ admits an eISS Lyapunov function of the form%
\begin{equation}\label{eq:Lyapunov-function-construction}
  V(x) = \sum_{i=1}^{\infty} \mu_i V_i(x_i),\quad V:X \rightarrow \R_+,%
\end{equation}
for some $\mu = (\mu_i)_{i\in\N}\in \ell^{\infty}$ satisfying $\underline{\mu} \leq \mu_i \leq \overline{\mu}$ with constants $\underline{\mu},\overline{\mu}>0$. In particular, the function $V$ has the following properties.%
\begin{enumerate}
\item[(a)] $V$ is continuous.%
\item[(b)] There is a $\lambda_\infty>0$ so that for all $x^0 \in X$ and $u \in \UC$%
\begin{equation*}
  \rmD^+ V_u(x^0) \leq -\lambda_{\infty} V(x^0) + \overline{\mu}\,\overline{\gamma}_u|u|_{q,\infty}^q.%
\end{equation*}
\item[(c)] For every $x \in X$ the following inequalities hold:%
\begin{equation}\label{eq:Coercivity-bound-for-V}
  \underline{\mu}\underline{\alpha}|x|_p^p \leq  V(x) \leq \overline{\mu}\,\overline{\alpha}|x|_p^p.%
\end{equation}
\end{enumerate}
In particular, $\Sigma$ is eISS.%
\end{theorem}

\begin{IEEEproof}
First, we prove the result for the case that there is a constant $\overline{\lambda}>0$ with%
\begin{equation}\label{eq:Upper-bound-for-lambda}
  \lambda_i \leq \overline{\lambda} \mbox{\quad for all\ } i\in\N.%
\end{equation}
Inequality \eqref{eq:Upper-bound-for-lambda} means that the decay rates of the eISS Lyapunov functions for all subsystems are uniformly bounded. Afterwards, we treat the general case.%

The proof proceeds in five steps.%

\emph{Step 1} (Definition and coercivity of $V$): 
First observe that the spectral radii of $\Psi = \Lambda^{-1}\Gamma:\ell^1 \rightarrow \ell^1$ and $\Theta = \Psi\trn:\ell^{\infty} \rightarrow \ell^{\infty}$ are the same, since the second operator is the adjoint of the first. Hence, Lemma \ref{lem_smallgain} yields a positive vector $\mu = (\mu_i)_{i\in\N} \in \ell^{\infty}$ whose components are uniformly bounded away from zero, and a constant $\lambda_{\infty}>0$ so that%
\begin{equation}\label{eq_muest}
  \frac{[\mu\trn(-\Lambda + \Gamma)]_i}{\mu_i} \leq -\lambda_{\infty} \quad \forall i \in \N.%
\end{equation}
To check that $V$ constructed as in \eqref{eq:Lyapunov-function-construction} is well-defined, note that for all $x \in X$ we have%
\begin{equation*}
  0 \leq V(x) \leq \sum_{i=1}^{\infty} \mu_i \overline{\alpha}_i |x_i|^p \leq \overline{\alpha}|\mu|_{\infty} |x|_p^p < \infty.%
\end{equation*}
This also shows the upper bound for \eqref{eq:Coercivity-bound-for-V}. The lower bound for \eqref{eq:Coercivity-bound-for-V} is obtained analogously, and thus inequality \eqref{eq_eiss_lyap_props-1} holds for $V$ (with $b=p$).%

\emph{Step 2} (Continuity of $V$): Fix a point $x = (x_i)_{i\in\N} \in X$ and some $\ep>0$. Choose $\delta_0,\ep'>0$ so that%
\begin{equation*}
  \overline{\mu}\,\overline{\alpha}2^{p-1}(\delta_0^p + \ep') \leq \frac{\ep}{4} \mbox{\quad and \quad} \ep' \leq \frac{\ep}{4\overline{\alpha}\,\overline{\mu}}.%
\end{equation*}
Subsequently, choose $N\in\N$ large enough such that $\sum_{i=N+1}^{\infty}|x_i|^p \leq \ep'$. Finally, choose $\delta \in (0,\delta_0]$ so that for every $y_i \in \R^{n_i}$ the following implication holds:%
\begin{equation*}
  |x_i - y_i| < \delta\ \Rightarrow\ |V_i(x_i) - V_i(y_i)| < \frac{\ep}{2N\overline{\mu}},\ i = 1,\ldots,N,%
\end{equation*}
where we use continuity of $V_i$ at $x_i$. Now let $y = (y_i)_{i\in\N} \in X$ be chosen so that $|x - y|_p < \delta$. In particular, this implies $|x_i - y_i| < \delta$ for $i=1,\ldots,N$. Then%
\begin{align*}
  |V(x) - V(y)| &\leq \sum_{i=1}^{\infty}\mu_i |V_i(x_i) - V_i(y_i)|\\
	&< \frac{\ep}{2} + \overline{\mu}\sum_{i=N+1}^{\infty} |V_i(x_i) - V_i(y_i)|.%
\end{align*}
The remainder sum can be estimated as%
\begin{align*}
  & \sum_{i=N+1}^{\infty} |V_i(x_i) - V_i(y_i)| \leq \overline{\alpha} \sum_{i=N+1}^{\infty}|x_i|^p + \overline{\alpha}\sum_{i=N+1}^{\infty} |y_i|^p\\
	&\leq \overline{\alpha}\ep' + \overline{\alpha}\sum_{i=N+1}^{\infty} 2^{p-1}\left(|y_i - x_i|^p + |x_i|^p\right)\\
	&\leq \overline{\alpha}\ep' + \overline{\alpha}2^{p-1}\left( \delta_0^p + \ep'\right) \leq \overline{\mu}^{-1}\frac{\ep}{2},%
\end{align*}
where we use $(a + b)^p \leq 2^{p-1}(a^p + b^p)$ for all $a,b\geq 0$, $p \geq 1$, which follows from the convexity of $a \mapsto a^p$. Altogether, one obtains $|V(x) - V(y)| < \ep$, showing that $V$ is continuous at $x$.%

\emph{Step 3} (Estimate of the orbital derivative): Fix an initial state $x^0 \in X$ and an external input $u\in\UC$. We write $\phi(t) = \phi(t,x^0,u)$, $\phi_i(t) = \pi_i\phi(t)$, $\bar{x}_i(t) = (\pi_j\phi(t))_{j\in I_i}$, where $\pi_i$ denotes the projection to the $i$-th component. Then for any $t>0$ (where $\phi(t)$ is defined), we obtain%
\begin{align*}
 \frac{1}{t}\big(V(\phi(t)) - V(x^0)\big) = \frac{1}{t}\sum_{i=1}^{\infty} \mu_i \big[V_i(\phi_i(t)) - V_i(\phi_i(0))\big].%
\end{align*}
Since the inequalities \eqref{eq_nablaviest} are valid for almost all positive times, the function on the right-hand side of \eqref{eq_nablaviest} is Lebesgue integrable, and since we assume that $\rmD_+(V_i \circ \phi_i)(t)<\infty$ for all $t$, we can proceed using the generalized fundamental theorem of calculus (see \cite[Thm.~9 and p.~42, Rmk.~5.c]{HaT06} or \cite[Thm.~7.3, p.~204]{Sak47}) to%
\begin{align*}
 \frac{1}{t}\big(V(\phi(t)) &- V(x^0)\big) \leq \frac{1}{t}\sum_{i=1}^{\infty}\int_0^t \mu_i \Bigl[-\lambda_i V_i(\phi_i(s)) \\
	&\qquad + \sum_{j\in I_i}\gamma_{ij}V_j(\phi_j(s)) + \gamma_{iu}|u_i(s)|^q\Bigr]\rmd s.%	
\end{align*}
We now apply the Fubini-Tonelli theorem in order to interchange the infinite sum and the integral (interpreting the sum as an integral associated with the counting measure on $\N$). To do this, it suffices to prove that the following integral is finite.%
\begin{equation*}
  \int_0^t \sum_{i=1}^{\infty} \Bigl|\mu_i \Bigl[-\lambda_i V_i(\phi_i(s)) + \sum_{j\in I_i}\gamma_{ij}V_j(\phi_j(s)) + \gamma_{iu}|u_i(s)|^q\Bigr]\Bigr| \rmd s.%
\end{equation*}
Using \eqref{eq_viest}, \eqref{eq_uniformity_alpha}, \eqref{eq_uniformity_gammaiu}, and \eqref{eq:Upper-bound-for-lambda}, we can upper bound the inner term by%
\begin{equation*}
   \overline{\mu}\Bigl[\overline{\lambda}\overline{\alpha}|\phi_i(s)|^p + \sum_{j \in I_i}\gamma_{ij}\overline{\alpha} |\phi_j(s)|^p 
	+ \overline{\gamma}_u|u_i(s)|^q\Bigr].%
\end{equation*}
By summing the three terms over $i$, one obtains%
\begin{align*}
  \overline{\lambda}\overline{\alpha} \sum_{i=1}^{\infty}|\phi_i(s)|^p &\leq c_1  |\phi(s)|_p^p,\\
	\sum_{i=1}^{\infty}\sum_{j \in I_i}\gamma_{ij}\overline{\alpha} |\phi_j(s)|^p &\leq c_2  \sum_{j=1}^{\infty} |\phi_j(s)|^p \sum_{i=1}^{\infty} \gamma_{ij} \leq c_3  |\phi(s)|_p^p,\\
  \overline{\gamma}_u\sum_{i=1}^{\infty}|u_i(s)|^q &= c_4  |u(s)|_{q}^{q}%
\end{align*}
for some constants $c_1,c_2,c_3,c_4>0$. In the inequality for the middle term, we use the boundedness assumption on the operator $\Gamma$. Hence,%
\begin{align*}
 & \int_0^t \sum_{i=1}^{\infty} \Bigl|\mu_i \Bigl[-\lambda_i V_i(\phi_i(s)) + \sum_{j\in I_i}\gamma_{ij}V_j(\phi_j(s))\\
&\qquad\qquad\qquad\qquad\qquad\qquad + \gamma_{iu}|u_i(s)|^q\Bigr]\Bigr| \rmd s\\
	&\quad \leq c  \int_0^t \left(|\phi(s)|_p^p + |u(s)|_{q}^q\right) \rmd s < \infty%
\end{align*}
for some constant $c>0$, where we use the fact that the integrand in the last term is essentially bounded ($s \mapsto |\phi(s)|_p^p$ is continuous and $s \mapsto |u(s)|_{q}^q$ is essentially bounded).%

Using the notation%
\begin{equation*}
  V_{\mathrm{vec}}(\phi(s)) := (V_1(\phi_1(s)),V_2(\phi_2(s)),\ldots)\trn%
\end{equation*}
and applying the Fubini-Tonelli theorem, we can then conclude that%
\begin{align*}
\frac{1}{t}\big(&V(\phi(t)) - V(x^0)\big)\\
&\leq \frac{1}{t}\int_0^t \sum_{i=1}^{\infty} \mu_i \Bigl[-\lambda_i V_i(\phi_i(s)) + \sum_{j\in I_i}\gamma_{ij}V_j(\phi_j(s)) \\
&\qquad\qquad\qquad\qquad\qquad\qquad\qquad\qquad + \gamma_{iu}|u_i(s)|^q\Bigr]\rmd s\\
	&= \frac{1}{t} \int_0^t \Bigl[ \mu\trn(-\Lambda + \Gamma)V_{\mathrm{vec}}(\phi(s)) + \sum_{i=1}^{\infty} \mu_i \gamma_{iu} |u_i(s)|^q\Bigr]\rmd s\\
	&\leq \frac{1}{t}\int_0^t \Bigl[-\lambda_{\infty} V(\phi(s)) + \overline{\mu}\,\overline{\gamma}_u|u|_{q,\infty}^q \Bigr]\rmd s\\
	&= \frac{1}{t}\int_0^t -\lambda_{\infty} V(\phi(s))\, \rmd s + \overline{\mu}\,\overline{\gamma}_u|u|_{q,\infty}^q,%
\end{align*}
where we use \eqref{eq_muest} to show the second inequality above. Since $s \mapsto V(\phi(s))$ is continuous, one obtains%
\begin{align*}
  \rmD^+ V_u(x^0) &= \limsup_{t \rightarrow 0+}\frac{1}{t}\left(V(\phi(t)) - V(x^0)\right) \\
  &\leq -\lambda_{\infty} V(x^0) + \overline{\mu}\,\overline{\gamma}_u|u|_{q,\infty}^q.%
\end{align*}
Hence, \eqref{eq_eiss_lyap_props-2} holds for $V$ with $\kappa = \lambda_{\infty}$ and $\gamma(r) = \overline{\mu}\,\overline{\gamma}_ur^q$.%

\emph{Step 4} (Proof of eISS): We showed that the properties (a)--(c) are satisfied for $V$. Thus, $V$ is an eISS Lyapunov function for $\Sigma$ and $\Sigma$ is eISS by Proposition \ref{prop_eiss}. Hence, for the case of uniformly upper-bounded $\lambda_i$ the theorem is proved.%

\emph{Step 5} (Unbounded decay rates $\lambda_i$): Assume that \eqref{eq:Upper-bound-for-lambda} does not hold for any $\overline{\lambda}$. Pick any $h>0$ and define the reduced decay rates $\lambda_i^h := \min\{\lambda_i,h\}$.
Thus, $\lambda_i^h \leq h$ for all $i\in\N$, which allows us to invoke the previous analysis.%

Indeed, as the inequalities \eqref{eq_nablaviest} hold with $\lambda_i$, they also hold with $\lambda_i^h$. Now let us define the modified operators $\Lambda^h$, $\Psi^h$ by%
\begin{equation*}
  \Lambda^h := \diag(\lambda_1^h,\lambda_2^h,\ldots),\qquad   \Psi^h := (\Lambda^h)^{-1}\Gamma.%
\end{equation*}
Considering $(\Lambda^h)^{-1}$ as an operator from $\ell^1$ to $\ell^1$, it is easy to see that $(\Lambda^h)^{-1} \to \Lambda^{-1}$ as $h\to\infty$. As $\Gamma$ is a bounded operator by assumption, it also holds that $\Psi^h \to \Psi$ as $h\to\infty$.%

By the small-gain condition, we have $r(\Psi)<1$. As the spectral radius is upper semicontinuous on the space of bounded operators on a Banach space (see e.g. \cite[Thm.~1.1(i)]{ADV08}), it holds that $r(\Psi^h) < 1$ for $h$ large enough. As the coefficients $\lambda_i^h$ are uniformly bounded, by feeding the operator $\Psi^h$ to Lemma \ref{lem_smallgain}, we obtain a vector $\mu = \mu(h)$ and a coefficient $\lambda_\infty = \lambda_\infty(h)$, so that (by the first four steps of this proof) \eqref{eq:Lyapunov-function-construction} is an eISS Lyapunov function for $\Sigma$ with decay rate $\lambda_{\infty}$.%
\end{IEEEproof}

\begin{remark}
Theorem~\ref{MT} provides a so-called dissipative form small-gain theorem. For large-but-finite networks, this form of SGCs has received considerable attention~\cite{Liu.2014,Ito.2013,Dashkovskiy.2011b} and has been applied to distributed control design~\cite{Liu.2014}, compositional construction of (in)finite-state abstractions~\cite{Pola.2016,Rungger.2016}, cyber-security of networked systems~\cite{Feng.2017}, and networked control systems with asynchronous communication~\cite{Heemels.2013}. Our result is a generalization of~\cite[Prop.~3.3]{Dashkovskiy.2011b} where the corresponding small-gain condition is a consequence of the Perron-Frobenius theorem; cf.~\cite[Lem.~3.1]{Dashkovskiy.2011b}. It basically relies on Lemma~\ref{lem_smallgain} which can be viewed as an infinite-dimensional extension of~\cite[Lem.~3.1]{Dashkovskiy.2011b}. \hfill$\diamond$%
\end{remark}

\begin{remark}\label{rem:Main-result} 
The parameters $\mu$ and $\lambda_\infty$, used for the construction of an eISS Lyapunov function $V$, are not uniquely determined. They depend, in particular, on the constant $\rho$ (in Lemma \ref{lem_smallgain}) and (in the case of unbounded $\lambda_i$) on the parameter $h$, introduced in Step 5 of the proof of Theorem~\ref{MT}.\hfill$\diamond$
\end{remark}

\subsection{Necessity of the required assumptions and tightness of the small-gain result}

The spectral small-gain condition $r(\Psi)<1$ cannot be removed or relaxed. This is already well-known for (nonlinear) planar systems; see~\cite[Sec.~1.5.4]{Mir12} for the tightness analysis of the small-gain condition.%

Assumptions \eqref{eq_viest} and \eqref{eq_uniformity_alpha} are necessary for the overall eISS Lyapunov function $V$ to be well-defined and coercive. If we remove the lower bound in \eqref{eq_viest} or \eqref{eq_uniformity_alpha}, we might still be able to prove ISS (though not eISS) of the interconnection by using results on non-coercive ISS Lyapunov functions, cf.~\cite[Thm.~2.18]{MiP19}, but additional assumptions on boundedness of reachable sets might be necessary.%

Without Assumption~\eqref{eq_uniformity_lambda} we do not have a uniform decay rate for the solutions of the subsystems, which prevents us from getting already asymptotic stability for the interconnection, even if the system is linear and all internal and external gains are zero. 
Consider, e.g., the infinite network
\begin{eqnarray*}
\dot{x}_i = -\frac{1}{i}x_i + u,\quad i\in\N,%
\end{eqnarray*}
with the input space $\UC := L^\infty(\R_+,\R)$, state space $X=\ell^p$ for any $p\in[1,\infty]$, and with Lyapunov functions $V_i(z)=z^2$ for all $i\in\N$ and $z\in\R$. With this choice of Lyapunov functions, all the assumptions which we impose for the small-gain theorem will be satisfied, except for \eqref{eq_uniformity_lambda}. At the same time, the network is not even exponentially stable in the absence of inputs. Furthermore, inputs of arbitrarily small magnitude may lead to unboundedness of trajectories, as mentioned, e.g., in~\cite[Sec.~6, p.~247]{MaP11}.%

Assumption~\eqref{eq_Gamma_bounded} is again crucial for the validity of the small-gain theorem, as shown by the following simple example%
\begin{eqnarray*}
  \dot{x}_i = -x_i + i u , \quad i \in\N,
\end{eqnarray*}
where we take $X=\ell^p$ for any $p\in[1,\infty]$. Choosing $V_i(z)=z^2$ for all $i\in\N$ and $z\in \R$, after some elementary manipulations one can see that all assumptions of the small-gain theorem are fulfilled, but the overall system is not ISS.%

\subsection{Verification of the spectral small-gain condition}\label{sec:Verification of the spectral small-gain condition}

The key requirement for the use of our small-gain theorem is a spectral small-gain condition $r(\Psi)<1$ for an operator $\Psi \in L(\ell^1)$, which is also equivalent to the condition $r(\Theta)<1$ for the bounded operator $\Theta=\Psi^*:\ell^\infty\to\ell^\infty$. Although there is no general, computationally efficient method to check this condition for any infinite-dimensional operator, there are many powerful equivalent criteria for this condition, which can be used for its verification.%

First of all, in view of Gelfand's formula, the spectral small-gain condition is equivalent to the existence of $k\in\N$ such that $\|\Psi^k\|<1$. This can be used in practice, especially for systems with a particular interconnection structure such as (quasi) spatially invariant systems. This is illustrated by several examples in Section~\ref{sec:Examples}.%

More generally, the condition $r(\Psi)<1$ is equivalent to the exponential stability of a discrete-time system, induced by the gain operator $\Psi$, for which an abundance of criteria exist, for a survey see \cite{GlM20}.

The previous criterion for $r(\Psi)<1$ (or $r(\Theta)<1$) can be used for arbitrary bounded operators $\Psi$, which do not have to be positive. In our case, $\Theta$ is a positive linear operator, acting on the space $\ell^{\infty}$, with the cone $\ell^\infty_+$, having a nonempty interior. With these strong additional properties, further powerful criteria for the spectral condition $r(\Theta)<1$ are available, which have been shown recently in \cite[Prop.~9, 10, 15]{MKG20} and \cite{GlM20}. 
More precisely, by \cite[Prop.~15]{MKG20}, \emph{the spectral small-gain condition $r(\Theta)<1$ is equivalent to the monotone bounded invertibility property of the operator $\id - \Theta$, which is equivalent due to \cite[Prop.~9, 10]{MKG20} to the so-called uniform small-gain condition of the operator $\Theta$}; we refer to \cite{MKG20} for the definitions of the corresponding notions. Finally, we note that $\ell^\infty$ with the cone $\ell^\infty_+$ is also a Banach lattice, which allows for several reformulations of the uniform small-gain conditions, see \cite[Rem.~5]{MKG20}.%

\section{Examples}\label{sec:Examples}

In this section we apply our results to three examples: linear spatially invariant systems, nonlinear spatially invariant systems with a nonlinearity satisfying the sector condition, and to a road traffic model. In all cases, we construct eISS Lyapunov functions with linear gains for all subsystems, and then apply our small-gain result to construct an exponential ISS Lyapunov function for the overall network.%

\subsection{A linear spatially invariant system}\label{exm:A-linear-spatially-invariant-system}

Consider an infinite network of systems $\Sigma_i$, given by%
\begin{align}
  \Sigma_i:\quad \dot x_i &= -b_{ii} x_i + b_{i(i-1)} x_{i-1} + b_{i(i+1)} x_{i+1}\nonumber\\
	&=: f_i(x_i,x_{i-1},x_{i+1}), \label{eq:spatially-invariant}
\end{align}
where $x_i \in \R$, $b_{ii} > 0$, $b_{i(i-1)},b_{i(i+1)} \in \R$ for each $i \in \N$ and $b_{10} = 0$.
We consider the standard Euclidean norm on $\R$ for each $i\in\N$ and assume that there is $\overline{b} > 0$ so that
\begin{equation*}
  \max\{b_{ii},|b_{i(i-1)}|,|b_{i(i+1)}|\} \leq \overline{b} \mbox{\quad for all\ } i \in \N.%
\end{equation*}
From Example \ref{ex_linear_general}, it immediately follows that the composite system is well-posed with $p = 2$.%

For each subsystem $\Sigma_i$, we choose the eISS Lyapunov function candidate $V_i(x_i) = \frac{1}{2}x_i^2$ satisfying~\eqref{eq_viest} and \eqref{eq_uniformity_alpha}.
Using Young's inequality, one can simply verify~\eqref{eq_nablaviest} as follows.
\begin{align*}
   \nabla V_i(x_i) \cdot& f_i (x_i,x_{i-1},x_{i+1})\\
	&= x_i(-b_{ii} x_i + b_{i(i-1)}x_{i-1} + b_{i(i+1)}x_{i+1}) \\
	&\leq -(b_{ii} - \ep_i - \delta_i) x_i^2 + \frac{b_{i(i-1)}^2}{4\ep_i} x_{i-1}^2 + \frac{b_{i(i+1)}^2}{4\delta_i}x_{i+1}^2 \\
	&= -2(b_{ii} - \ep_i - \delta_i)V_i(x_i) + \frac{b_{i(i-1)}^2}{2\ep_i} V_{i-1}(x_{i-1})\\
	&\qquad + \frac{b_{i(i+1)}^2}{2\delta_i}V_{i+1}(x_{i+1}),%
\end{align*}
for appropriate choices of $\ep_i,\delta_i > 0$. Hence, we can choose
\begin{align*}
  \lambda_i &:= 2(b_{ii} - \ep_i - \delta_i),~\gamma_{i(i-1)} := \frac{b_{i(i-1)}^2}{2\ep_i},~\gamma_{i(i+1)} := \frac{b_{i(i+1)}^2}{2\delta_i},%
\end{align*}
and assume that $\ep_i,\delta_i$ are such that \eqref{eq_uniformity_lambda} and \eqref{eq_Gamma_bounded} are satisfied.
It follows that the infinite matrix $\Psi$ has the form%
{\small
\begin{equation}\label{eq:Psi-example}
  \Psi {=} \Lambda^{-1}\Gamma {=} \left(\begin{array}{ccccccc}
0 & \psi_{12} & 0 & 0 & 0 & 0 & \dots \\
\psi_{21} & 0 & \psi_{23} & 0 & 0 & 0 & \dots \\
0 & \psi_{32} & 0 & \psi_{34} & 0 & 0 & \dots \\
0 & 0 & \psi_{43} & 0 & \psi_{45} &  0 & \dots \\
\vdots & \ddots & \ddots & \ddots & \ddots & \ddots & \ddots
\end{array}\right),%
\end{equation}}
where $\psi_{ij} = \gamma_{ij}/\lambda_i$.

In view of Theorem~\ref{MT}, the following set of sufficient conditions guarantees that the interconnection defined above is eISS.%:%
\begin{itemize}
\item $\max\{b_{ii},|b_{i(i-1)}|,|b_{i(i+1)}|\} \leq \overline{b}$ for all $i\in\N$ with a constant $\overline{b}>0$ (for well-posedness).%
\item Constants $\ep_i,\delta_i>0$ are chosen such that%
\begin{itemize}
	\item Assumption (ii) in Theorem~\ref{MT} holds, i.e.:\linebreak
	 $\exists \underline{\lambda}>0$ with $\underline{\lambda} \leq 2(b_{ii} - \ep_i - \delta_i)$ for all $i\in\N$. 
	\item Assumption (iii) in Theorem~\ref{MT} holds, i.e.:\linebreak
	$\exists \overline{\gamma}<\infty$: $\frac{b_{i(i-1)}^2}{2\ep_i}\leq \overline{\gamma}$ and $\frac{b_{i(i+1)}^2}{2\delta_i} \leq \overline{\gamma}$ for all $i\in\N$.
	\item The small-gain condition $r(\Psi)<1$ holds. A sufficient condition for this would be
\begin{equation*}
  r(\Psi) \leq \|\Psi\| = \sup_{j\in\N} \sum_{i=1}^{\infty} \psi_{ij} < 1,%
\end{equation*}	
which is equivalent to existence of $\tilde{\psi}<1$ such that for all $i\in\N$
\begin{align*}
  \psi_{i(i-1)}+\psi_{i(i+1)} &= \frac{b_{i(i-1)}^2}{4\ep_i(b_{ii} {-} \ep_i {-} \delta_i)}\\
	&+ \frac{b_{i(i+1)}^2}{4\delta_i(b_{ii} {-} \ep_i {-} \delta_i)}\leq \tilde{\psi}.
\end{align*}
\end{itemize}
\end{itemize}

\begin{remark}\label{rem:Example-A-} 
As we argued in Remark~\ref{rem:Choice-of-a-state space}, the choice of the ``right'' Lyapunov functions depends on the physical meaning of the variables, and thus quadratic Lyapunov functions, and the corresponding state space $X = \ell^2(\N,(n_i))$ may not be physically appropriate for some applications. However, there are other natural options for Lyapunov functions for the subsystems $\Sigma_i$, for example $W_i(x_i):=|x_i|$, which would lead to other values of the gains, and to another expression for the small-gain condition. \hfill$\diamond$%
\end{remark}

\subsection{A nonlinear multidimensional spatially invariant system}%

Here, we analyze a class of nonlinear control systems which widely appeared in many applications, including neural networks, analysis and design of optimization algorithms, Lur'e problem, and so on (see \cite{fetzer} and the references therein).%

Consider an infinite network whose subsystems are described by%
\begin{equation*}
  \Sigma_i:\quad \dot x_i = A_{i} x_i + E_i\varphi_i(G_ix_i) + B_iu_i+D_i\bar{x}_i,%
\end{equation*}
where $A_i\in \R^{n_i\times n_i}$, $E_i\in \R^{n_i}$, $G^{\trn}_i\in \R^{n_i}$, $B_i\in \R^{n_i\times m_i}$, $D_i\in \R^{n_i\times N_i}$ with $N_i=\sum_{j\in I_i}n_j$ and $I_1 = \{i+1\}$, $I_i = \{i-1,i+1\}$ for all $i \geq 2$.%

We consider the standard Euclidean norm on each $\R^{n_i}$, $\R^{m_i}$ and $\R^{N_i}$, and assume that $A_i,E_i,G_i,B_i,D_i$ are uniformly bounded for all $i\in \N$. That is, $\Vert A_i\Vert\leq a, \Vert E_i\Vert\leq e,\Vert G_i\Vert\leq g ,\Vert B_i\Vert\leq b ,\Vert D_i\Vert\leq d$. Additionally, we assume that the nonlinear functions $\varphi_i:\R\rightarrow \R$ satisfy%
\begin{align}\label{sect}
  \big(\varphi_i(G_ix_i)-r_iG_ix_i\big)\big(\varphi_i(G_ix_i)-l_iG_ix_i\big) \leq 0%
\end{align} 
for all $x_i\in\R^{n_i}$ with $r_i>l_i$, $l_i,r_i\in\R$. Moreover, we assume that the nonlinear functions $\varphi_i:\R\rightarrow \R$ have some regularity properties such that the interconnected system $\Sigma$ with state space $X:=\ell^2(\N,(n_i))$ and input space $U:=\ell^2(\N,(m_i))$ is well-posed.%

Now let for all $i\in\N$ the function $V_i$ be defined as $V_i(x_i) := x_i\!{\trn} M_ix_i$, where $M_i\in\R^{n_i\times n_i}$ is a symmetric and positive definite matrix with $\Vert M_i\Vert\leq m$ and $0<\underline{m}\leq\lambda_{\min}(M_i)\leq\lambda_{\max}(M_i)\leq \overline{m}<\infty$, where $\lambda_{\min}(\cdot)$ and $\lambda_{\max}(\cdot)$ denote the smallest and largest eigenvalues, respectively.%

Assume that for all $i\in\N$, $x_i \in \R^{n_i}$, and $\varphi_i:\R\rightarrow \R$ satisfying \eqref{sect}, the inequality%
\begin{equation}\label{lmi}
  2x_i\trn M_i(A_i x_i + E_i\varphi_i(G_i x_i)) \leq -\kappa_i x_i\trn M_i x_i%
\end{equation} 
holds for some $\kappa_i$ with $0 < \underline{\kappa} \leq \kappa_i$ for some $\underline{\kappa}$.%

\begin{remark}
Note that inequality \eqref{lmi} is equivalent to%
\begin{align*}
\begin{bmatrix}
x_i\\
\varphi_i(G_ix_i) 
\end{bmatrix}\trn\begin{bmatrix}
A\trn_i M_i\!\!+\!\!M_iA_i\!\!+\!\!\kappa_i M_i& M_iE_i\\
E\trn_iM_i& 
\textbf{0}
\end{bmatrix}\begin{bmatrix}
x_i\\
\varphi_i(G_ix_i)
\end{bmatrix}&\! \leq 0%
\end{align*}
for all $i\in\N$ and $x_i \in \R^{n_i}$, where $\textbf{0}$ is a zero matrix of appropriate dimensions. Now note that inequality \eqref{sect} is equivalent to%
\begin{align}\label{sectlmi}
\begin{bmatrix}
x_i\\
\varphi_i(G_ix_i) 
\end{bmatrix}\trn\begin{bmatrix}
r_il_iG\trn_i G_i& \!-\frac{r_i+l_i}{2}G\trn_i\\
-\frac{r_i+l_i}{2}G_i& 
\!1
\end{bmatrix}\begin{bmatrix}
x_i\\
\varphi_i(G_ix_i)
\end{bmatrix}&\!\! \leq \!0.
\end{align}
Hence, by using \eqref{sectlmi} and the S-procedure \cite{slemma}, a sufficient condition for the validity of \eqref{lmi} is the validity of  the matrix inequality%
\begin{align*}
\begin{bmatrix}
A\trn_i M_i\!+\!M_iA_i\!+\!\kappa_i M_i-r_il_iG\trn_i G_i& M_iE_i\!+\!\tau_i\frac{r_i+l_i}{2}G\trn_i\\
E\trn_iM_i\!+\!\tau_i\frac{r_i+l_i}{2}G_i& -\tau_i
\end{bmatrix}&\!\!\preceq\!0%
\end{align*}
for some $\tau_i\in\R_+$. \hfill$\diamond$%
\end{remark}

Then%
\begin{equation*}
  \lambda_{\min}(M_i)|x_i|^2\leq V_i(x_i)\leq\lambda_{\max}(M_i)|x_i|^2
\end{equation*} 
and%
\begin{align*}
\nabla V_i&(x_i) \cdot f_i(x_i,\bar{x}_i,u_i)\\
&= 2x_i\trn M_i(A_i x_i +E_i\varphi_i(G_ix_i) + B_iu_i + D_i\bar{x}_i)\\
&= 2x_i\trn M_i \big(A_ix_i + E_i\varphi_i(G_ix_i)\big)\\
&\qquad\qquad\qquad\qquad\qquad+ 2x_i\trn M_iB_iu_i + 2x_i\trn M_iD_i\bar{x}_i.%
\end{align*}
Using Cauchy-Schwarz and Young's inequalities, respectively, we obtain (for any $\ep_i>0$)%
\begin{align*}
2x_i\trn M_iB_iu_i &= 2x_i\trn \sqrt{M_i}\sqrt{M_i} B_iu_i \\
&\leq 2|\sqrt{M_i}x_i|\cdot|\sqrt{M_i}B_iu_i| \\
&\leq 2|\sqrt{M_i}x_i|\cdot\|\sqrt{M_i}B_i\| |u_i| \\
&\leq \ep_i|\sqrt{M_i}x_i|^2 + \frac{\|\sqrt{M_i}B_i\|^2|u_i|^2}{\ep_i} \\
&= \ep_i x_i\trn M_i x_i + \frac{\|\sqrt{M_i}B_i\|^2|u_i|^2}{\ep_i},%
\end{align*}
and analogously,%
\begin{equation}
\label{eq:Estimate-xbar}
2x_i\trn M_iD_i\bar{x}_i \leq \ep_i x_i\trn M_i x_i + \frac{\|\sqrt{M_i}D_i\|^2|\bar{x}_i|^2}{\ep_i}.%
\end{equation}
Therefore, we have%
\begin{align*}
& \nabla V_i(x_i) \cdot f_i(x_i,\bar{x}_i,u_i) \leq -(\kappa_i - 2\ep_i) x_i\trn M_i x_i\nonumber\\
&\qquad + \frac{\|\sqrt{M_i}B_i\|^2|u_i|^2}{\ep_i} + \frac{\|\sqrt{M_i}D_i\|^2}{\ep_i}(|x_{i-1}|^2 + |x_{i+1}|^2)\nonumber\\
&\leq -(\kappa_i - 2\ep_i) V_i(x_i) + \frac{\|\sqrt{M_i}B_i\|^2}{\ep_i}|u_i|^2 \nonumber\\
&\qquad + \frac{\|\sqrt{M_i}D_i\|^2}{\ep_i}\Bigl(\frac{V_{i-1}(x_{i-1})}{\lambda_{\min}(M_{i-1})} + \frac{V_{i+1}(x_{i+1})}{\lambda_{\min}(M_{i+1})}\Bigr).%
\end{align*}
Hence, the function $V_i(x_i) = x_i\!{\trn}M_ix_i$ is an eISS Lyapunov function for the subsystem $\Sigma_i$ satisfying \eqref{eq_viest} and \eqref{eq_nablaviest} with%
\begin{align*}
\underline{\alpha}_i &:= \lambda_{\min}(M_i),~~\overline{\alpha}_i := \lambda_{\max}(M_i),~~ \lambda_i := \kappa_i - 2\ep_i,\\
\gamma_{ij} &:= \frac{\Vert\sqrt{M_i}D_i\Vert^2}{\lambda_{\min}(M_j)\ep_i},~~\gamma_{iu} := \frac{\Vert \sqrt{M_i}B_i\Vert^2}{\ep_i}.%
\end{align*}
With $\underline{\alpha} := \underline{m}$ and $\overline{\alpha} := \overline{m}$, \eqref{eq_uniformity_alpha} is satisfied. With a uniformity condition on $\ep_i$, say $0 < \underline{\ep} \leq \ep_i \leq \overline{\ep}<\infty$ so that $\underline{\kappa} - 2\overline{\ep} > 0$, we see that \eqref{eq_uniformity_lambda} also holds with $\underline{\lambda} := \underline{\kappa} - 2\overline{\ep}$. Finally, we have%
\begin{equation*}
0 < \gamma_{ij} \leq \frac{md^2}{\underline{m}\,\underline{\ep}} =: \overline{\gamma} < \infty,%
\end{equation*}
showing that \eqref{eq_Gamma_bounded} is satisfied by Lemma \ref{lem_Gamma_bounded}, and%
\begin{equation*}
\gamma_{iu} \leq \frac{mb^2}{\underline{\ep}} =: \overline{\gamma}_u \mbox{\quad for all\ } i \in \N,%
\end{equation*}
which implies \eqref{eq_uniformity_gammaiu}.
Clearly, the infinite matrix $\Psi := \Lambda^{-1}\Gamma$, for $\Lambda $ and $\Gamma$ as in~\eqref{eq:operator-A}, has the same form as the one in~\eqref{eq:Psi-example}.%

In that way, with the same arguments as in Section~\ref{exm:A-linear-spatially-invariant-system}, one can conclude that any choice of the numbers $\ep_i$ such that $\frac{md^2}{(\underline{\kappa}-2\overline{\ep})\underline{m}\,\underline{\ep}} < \frac{1}{2}$ for all $i\in\N$ leads to $r(\Psi) < 1$.%

Hence, by Theorem \ref{MT} there exists $\mu = (\mu_i)_{i\in\N}\in \ell^{\infty}$ satisfying $0 < \underline{\mu} \leq \mu_i \leq \overline{\mu} < \infty$ with constants $\underline{\mu},\overline{\mu}$ such that function $V(x) = \sum_{i=1}^{\infty} \mu_i x_i\!{\trn} M_ix_i$ is an eISS Lyapunov function for the interconnected system $\Sigma$.%

We note that our example here covers the most common mathematical representation of the tracking error in a vehicle platooning application. The choice of%
\begin{equation*}
A = \left[\begin{array}{cc}
0 & 1 \\
-k_0 & -b_0
\end{array}\right],\
B_i = \left[\begin{array}{c}
0 \\
1
\end{array}\right],\
D = \left[\begin{array}{cc}
0 & 0 \\
k_0 & b_0
\end{array}\right],%
\end{equation*}
where $k_0,b_0 > 0$, $G(x_i) = 0$, gives the closed-loop dynamics for the symmetric bidirectional architecture in a vehicle platooning problem, where the control action of each vehicle depends on the relative position and velocity measurements from its immediate front and back neighbors~\cite{Hao.2013,Jovanovic.2005b,BPD02}.%

\begin{remark}\label{rem:Alternative-norms-for-xbar}
In this example, we have used the Euclidean norm for the space $\R^{N_i}$ of $\bar{x}_i$-vectors. However, one could utilize in computations also more specific norms, based on our choice of the Lyapunov functions $V_i$, which lead to more precise bounds on the gains.%

For example, note that $\sqrt{V_{i-1}(\cdot)}$ is a norm on $\R^{n_i}$, and thus one can define a norm on $\R^{N_i}$ by%
\begin{equation*}
  \bar{x}_i \mapsto \normt{\bar{x}_i}:= \sqrt{V_{i-1}(x_{i-1}) + V_{i+1}(x_{i+1})}.%
\end{equation*}
Now instead of \eqref{eq:Estimate-xbar}, we could obtain the estimate%
\begin{equation*}
  2x_i\trn M_iD_i\bar{x}_i \leq \ep_i x_i\trn M_i x_i + \frac{\|\sqrt{M_i}D_i\|_{\normt{\cdot},2}^2\normt{\bar{x}_i}^2}{\ep_i},%
\end{equation*}
where the double index on the left-hand side indicates that we consider the operator norm induced by the norm $\normt{\cdot}$ on the domain and by Euclidean norm on the codomain of the corresponding operator. Proceeding further, we can again verify that the function $V_i(x_i) = x_i\!{\trn}M_ix_i$ is an eISS Lyapunov function for the subsystem $\Sigma_i$, and at the same time avoid a rather rough estimate of $|x_i|^2$ by $\frac{V_{i-1}(x_{i-1})}{\lambda_{\min}(M_{i-1})}$. \hfill$\diamond$%
\end{remark}

\subsection{A road traffic model}

\begin{figure}
	\vspace*{-0.0cm}
	\begin{center}
			\includegraphics[height=5.5cm]{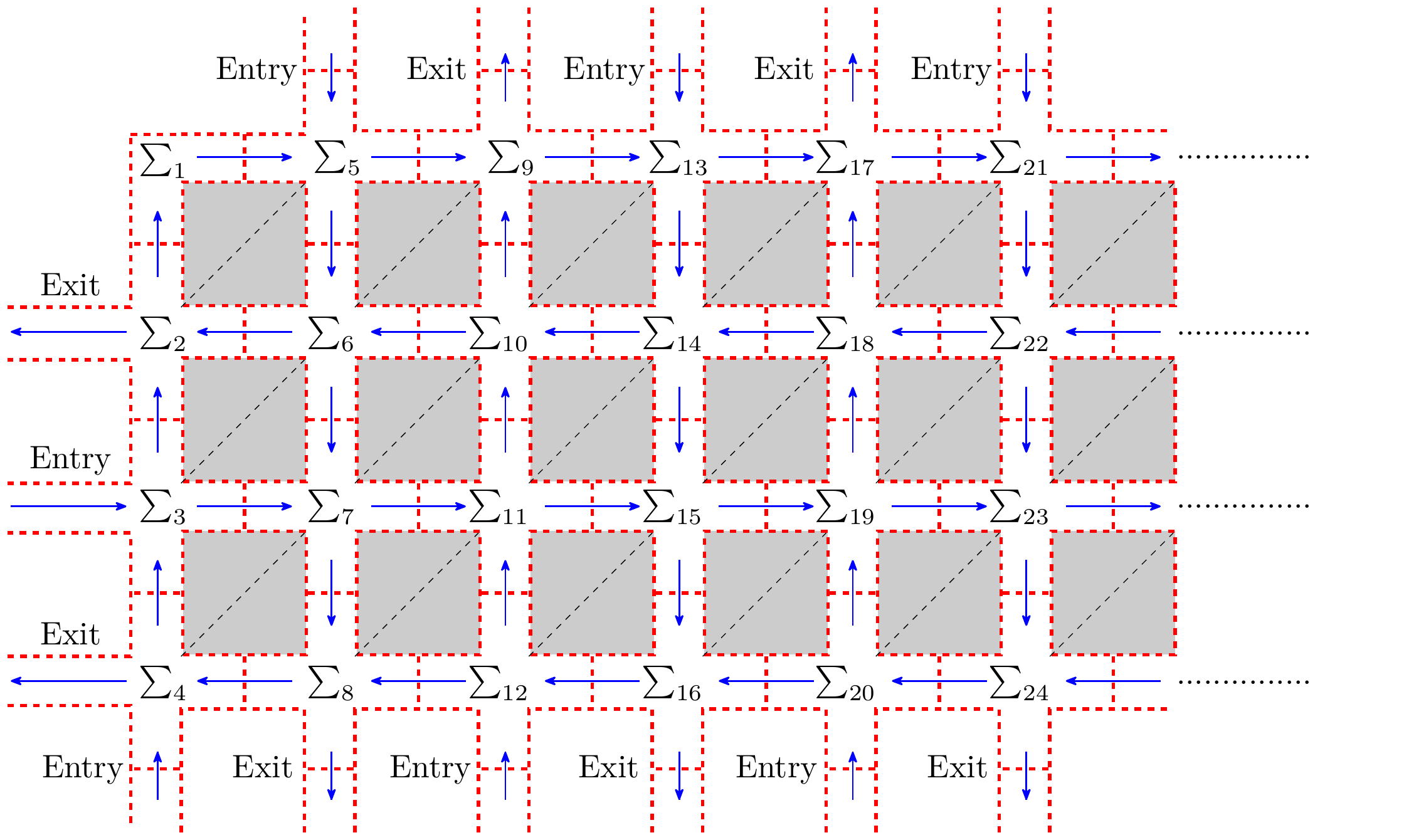}
				\vspace*{-0.6cm}
		\caption{Model of a road traffic network composed of infinitely many subsystems.}
		\label{allr}
	\end{center}
	\vspace*{-0.2cm}
\end{figure}

In this example, we apply our approach to a variant of the road traffic model from \cite{Kibangou}.
We consider a traffic network divided into infinitely many cells, indexed by $i\in\N$.
Each cell $i$ represents a subsystem $\Sigma_i$ described by a differential equation of the following form%
\begin{align}\label{subsys}
  \Sigma_i: \dot x_i = -\Bigl(\frac{v_i}{l_i} + e_i\Bigr)x_i + D_i\bar{x}_i + B_iu_i,\quad x_i,u_i\in \R,%
\end{align}
with the following structure%
\begin{enumerate}
\item[$-$] $e_i=0,D_i=c\frac{v_{i+1}}{l_{i+1}},\bar{x}_i=x_{i+1},B_i=0$ if $i\in S_1:=\{1,3\}$;
\item [$-$]$e_i=0,D_i=c\frac{v_{{i+4}}}{l_{i+4}},\bar{x}_i=x_{i+4},B_i=r>0$ if $i\in S_2:=\{4+8k : k\in\N\cup\{0\}\}$;
\item [$-$]$e_i=0,D_i=c\frac{v_{{i-4}}}{l_{i-4}},\bar{x}_i=x_{i-4},B_i=\frac{r}{2}$ if $i\in S_3:=\{5+8k :k\in\N\cup\{0\}\}$;
\item [$-$]$e_i=0,D_i=c(\frac{v_{i-1}}{l_{i-1}},\frac{v_{i+4}}{l_{i+4}})\trn,\bar{x}_i=(x_{i-1},x_{i+4}),B_i=0$ if $i\in S_4:=\{6+8k : k\in\N\cup\{0\}\}$;
\item [$-$] $e_i=e\in(0,1),D_i=c(\frac{v_{i-4}}{l_{i-4}},\frac{v_{i+1}}{l_{i+1}})\trn,\bar{x}_i=(x_{i-4},x_{i+1}),B_i=0$ if $i\in S_5:=\{9+8k : k\in\N\cup\{0\}\}$;
\item [$-$] $e_i=0,D_i=c(\frac{v_{i+1}}{l_{i+1}},\frac{v_{i+4}}{l_{i+4}})\trn,\bar{x}_i=(x_{i+1},x_{i+4}),B_i=0$ if $i\in S_6:=\{2+8k : k\in\N\cup\{0\}\}$;
\item [$-$] $e_i=0,D_i=c(\frac{v_{i-4}}{l_{i-4}},\frac{v_{i-1}}{l_{i-1}})\trn,\bar{x}_i=(x_{i-4},x_{i-1}),B_i=0$ if $i\in S_7:=\{7+8k : k\in\N\cup\{0\}\}$;
\item [$-$] $e_i=2e,D_i=c(\frac{v_{i-1}}{l_{i-1}},\frac{v_{i+4}}{l_{i+4}})\trn,\bar{x}_i=(x_{i-1},x_{i+4}),B_i=0$ if $i\in S_8:=\{8+8k : k\in\N\cup\{0\}\}$;
\item [$-$] $e_i=0,D_i=c(\frac{v_{i-4}}{l_{i-4}},\frac{v_{i+1}}{l_{i+1}})\trn,\bar{x}_i=(x_{i-4},x_{i+1}),B_i=0$ if $i\in S_9:=\{11+8 k : k\in\N\cup\{0\}\}$;
\end{enumerate}
where, for all $i\in\N$, $0~<\underline{v}\leq~ v_i\leq \overline{v}$, $0<\underline{l}\leq l_i\leq \overline{l}$, and $c\in(0,0.5)$. In \eqref{subsys}, $l_i$ is the length of a cell in kilometers (km), and $v_i$ is the flow speed of the vehicles in kilometers per hour (km/h).
The state of each subsystem $\Sigma_i$, i.e. $x_i$, is the density of traffic, given in vehicles per cell, for each cell $i$ of the road.
The scalars $B_i$ represent the number of vehicles that can enter the cells through entries which are controlled by $u_i$.
In particular, $u_i=1$ means green light and $u_i=0$ means red light. Moreover, the constants $e_i$ represent the percentage of vehicles that leave the cells using available exits. The overall system and subsystems are illustrated by Figure~\ref{allr}.%

%\begin{figure}
	%\begin{center}
		%\includegraphics[height=5.5cm]{rte-eps-converted-to}
		%\caption{Model of a road traffic network composed of infinitely many subsystems.}
		%\label{allr}
	%\end{center}
	%\vspace*{-0.7cm}
%\end{figure}
Clearly, the interconnected system $\Sigma$ with state space $X:=\ell^2(\N,(n_i))$ and input space $U:=\ell^{2}(\N,(m_i))$ is well-posed (cf.~Example \ref{ex_linear_general}). The choice $U:=\ell^{2}(\N,(m_i))$ means that only a finitely many of the traffic lights can be green at each $t \geq 0$.%

Furthermore, each subsystem $\Sigma_i$ admits an eISS Lyapunov function of the form $V_i(x_i) = \frac{1}{2}x_i^2$. The function $V_i$ satisfies \eqref{eq_viest} and \eqref{eq_nablaviest} for all $i\in\N$ with $\underline{\alpha}_i = \overline{\alpha}_i = \frac{1}{2}$, $\lambda_i = 2(\frac{v_i}{l_i}+e_i-2\ep_i) $, $\gamma_{ij} = \frac{\Vert D_i\Vert^2}{2\ep_i}$ for all $j\in I_i$, $\gamma_{iu} = \frac{B_i^2}{2\ep_i}$, for an appropriate choice of $0<\underline{\ep} \leq \ep_i \leq \overline{\ep}$ such that $ 0 < \underline{\lambda} := 2\bigl(\underline{v}/\overline{l} - 2\overline{\ep}\bigr) \leq \lambda_i$.
In that way, one can readily observe that%
\begin{align*}
 	0 < \gamma_{ij} \leq \frac{(c\overline{v})^2}{\underline{\ep}\,\underline{l}^2} =: \overline{\gamma} < \infty,\quad0 < \gamma_{iu} \leq \frac{r^2}{2\underline{\ep}} =: \overline{\gamma}_u < \infty.%
\end{align*}

Additionally, the infinite matrix $\Psi := \Lambda^{-1}\Gamma = (\psi_{ij})_{i,j\in\N} = (\gamma_{ij}/\lambda_i)_{i,j\in\N}$, for $\Lambda $ and $\Gamma$ defined in \eqref{eq:operator-A}, has the following structure.%:%
\begin{enumerate}
	\item[$-$] $i \in S_1$ $\Rightarrow$ ($\gamma_{ij}\neq0$ $\Leftrightarrow$ $j=i+1$);%
	\item[$-$] $i \in S_2$ $\Rightarrow$ ($\gamma_{ij}\neq0$ $\Leftrightarrow$ $j=i+4$);%
	\item [$-$] $i \in S_3$ $\Rightarrow$ ($\gamma_{ij}\neq0$ $\Leftrightarrow$ $j=i-4$);%
	\item[$-$] $i \in S_4$ $\Rightarrow $ ($\gamma_{ij}\neq0$ $\Leftrightarrow$ $j\in \{i-1,i+4\}$);%
	\item [$-$] $i \in S_5$ $\Rightarrow$ ($\gamma_{ij}\neq0$ $\Leftrightarrow$ $j\in \{i-4,i+1\}$);%
	\item [$-$] $i \in S_6$ $\Rightarrow$ ($\gamma_{ij}\neq0$ $\Leftrightarrow$ $j\in \{i+1,i+4\}$);%
	\item [$-$] $i \in S_7$ $\Rightarrow$ ($\gamma_{ij}\neq0$ $\Leftrightarrow$ $j\in \{i-4,i-1\}$);%
	\item [$-$] $i \in S_8$ $\Rightarrow$ ($\gamma_{ij}\neq0$ $\Leftrightarrow$ $j\in \{i-1,i+4\}$);%
	\item [$-$] $i \in S_9$ $\Rightarrow$ ($\gamma_{ij}\neq0$ $\Leftrightarrow$ $j\in \{i-4,i+1\}$).%
\end{enumerate}
The spectral radius $r(\Psi)$ can be estimated by%
\begin{align*}
  r(\Psi) \leq \|\Psi\| = \sup_{j\in\N} \sum_{i=1}^{\infty} \psi_{ij} \leq 2 \frac{\overline{\gamma}}{\underline{\lambda}}.%
\end{align*}

Hence, any choice of the constants $\ep_i$ such that%
\begin{align*}
({2(c\overline{v})^2}/{\underline{\ep}~\underline{l}^2})/(({\underline{v}}/{\overline{l}})-2\overline{\ep}) < 1,%
\end{align*}
for all $i\in\N$, leads to $r(\Psi) < 1$.%

Hence, by Theorem \ref{MT} there exists $\mu = (\mu_i)_{i\in\N}\in \ell^{\infty}$ satisfying $\underline{\mu} \leq \mu_i \leq \overline{\mu}$ with constants $\underline{\mu},\overline{\mu}>0$ such that the function $V(x) = \frac{1}{2}\sum_{i=1}^{\infty} \mu_i x_i^2$ is an eISS Lyapunov function for the interconnected system $\Sigma$.%

\section{Conclusions}\label{sec:Conclusions}

In this paper, we developed sufficient small-gain type conditions for showing exponential ISS of networks consisting of  countably infinite numbers of exponentially ISS subsystems. Our main mathematical tool is the theory of positive linear operators in ordered Banach spaces. The proposed small-gain condition, expressed in terms of the spectral radius of the resulting gain operator, has powerful criteria in terms of monotone bounded invertibility, uniform small-gain conditions and other stability properties, as shown in \cite{MKG20, GlM20}. We applied our results to some classes of linear and nonlinear systems.%

Our results can be extended in several directions. A challenging open question is whether similar conditions can be derived for (generally, non-exponential) input-to-state stability of countable interconnections of merely ISS subsystems.%

Nonlinear small-gain theorems in trajectory formulation for infinite networks, obtained in \cite{MKG20}, as well as novel relations between the monotone limit property, spectral and uniform small-gain conditions \cite{MKG20}, provide a strong basis for tackling this problem.%

Another direction is to use our results for the development of scale-free distributed/decentralized control design, which is now under investigation. In the spirit of~\cite{Noroozi.2018a}, our future work also investigates the necessity of such small-gain conditions.%

Furthermore, it is interesting to prove the small-gain theorem for networks, modeled on $\ell^\infty$ (as opposed to the $\ell^p$ space for finite $p$). The proof technique, used in this paper is not easily adaptable to the $\ell^\infty$-case, and probably other proof techniques are required, e.g., max-type ISS Lyapunov functions for interconnections may be more appropriate in this case, see \cite{DMS19a}.%

\section*{Acknowledgment}

The authors would like to thank Fabian Wirth for insightful discussions on small-gain theorems and Jochen Gl\"{u}ck for helping us to simplify the proof of Lemma \ref{lem_smallgain}.%

\bibliographystyle{IEEEtran}
\bibliography{references}

\appendices

\section{Positive operators}\label{sec_positive_operators}

In this section, we recall some results about positive operators on ordered Banach spaces. We start with some elementary facts about bounded operators. A general reference is \cite{DunSch57}.%

Let $X$ be a real Banach space with norm $|\cdot|$ and $T:X \rightarrow X$ a bounded linear operator. Recall that the complexification of $X$ is the complex Banach space $X_{\C} = \{x + iy : x,y\in X\}$ equipped with the norm $|x+iy| := \sup_{t\in[0,2\pi]}|(\cos t)x + (\sin t)y|$. The complexification of $T$ is the bounded operator $T_{\C}(x+iy) := Tx + iTy$, $T_{\C}:X_{\C} \rightarrow X_{\C}$. The resolvent of $T$ is the function $R(\lambda,T) := (\lambda I - T_{\C})^{-1}$, defined for all $\lambda \in \C$ so that the inverse exists and is a bounded operator. The \emph{resolvent set} of $T$ is $\rho(T) = \{\lambda \in \C : R(\lambda,T) \mbox{ exists and is bounded}\}$ and the \emph{spectrum} of $T$ is $\sigma(T) = \C \backslash \rho(T)$, which is a nonempty compact set. The \emph{spectral radius} of $T$ is defined as%
\begin{equation*}
  r(T) := \max\{ |\lambda| : \lambda \in \sigma(T) \}.%
\end{equation*}
A way to compute $r(T)$ is provided by Gelfand's formula \cite{DunSch57}:%
\begin{equation}\label{eq_gelfand}
  r(T) = \lim_{n \rightarrow \infty}\|T^n\|^{1/n} = \inf_{n\in\N}\|T^n\|^{1/n},%
\end{equation}
where $\|\cdot\|$ denotes the operator norm induced by the norm on $X$. Let $X^*$ denote the topological dual space of $X$, i.e., the Banach space of all bounded linear functionals $x^*:X \rightarrow \R$, equipped with the operator norm $|x^*| = \sup_{|x|=1}|x^*(x)|$. The \emph{adjoint operator} of $T$ is defined as $(T^*x^*)(x) := x^*(Tx)$ for all $x^* \in X^*$ and $x \in X$. The adjoint operator satisfies $\sigma(T^*) = \sigma(T)$ and $\|T^*\| = \|T\|$.%

A nonempty subset $K \subset X$ is called a \emph{cone} if it is closed\footnote{Sometimes the closedness is not part of the definition of a cone, and cones satisfying this assumption are called closed cones.} and convex and satisfies the following properties:%
\begin{itemize}
\item If $x \in K$ and $\lambda \geq 0$, then $\lambda x \in K$.%
\item If $x,-x\in K$, then $x = 0$.%
\end{itemize}
In particular, the former of these properties together with the convexity implies that for any $x,y \in K$ and $\lambda,\mu \geq 0$ also $\lambda x + \mu y \in K$.%

The specification of a cone in $X$ defines a partial order $\geq$ for all $x,y\in X$ by $x \geq y$ if and only if $x - y \in K$. The pair $(X,K)$ is thus called an \emph{ordered Banach space}.%

Once a cone $K$ has been specified, a bounded linear operator $T:X \rightarrow X$ is called \emph{positive} if $T(K) \subset K$. In this case, we write $T>0$ if $T \neq 0$. The positivity of an operator $T$ can also be expressed by the implication%
\begin{equation*}
  x \geq y \quad\Rightarrow\quad Tx \geq Ty,\quad x,y \in X.%
\end{equation*}

If $\inner(K)\neq\emptyset$, we write $x\ll y$ if $y-x \in\inner(K)$.%

\begin{lemma}\label{lem:rT_bound-new}\footnote{Lemma~\ref{lem:rT_bound-new} as well as its proof was proposed to us by Jochen Gl\"uck, see also \cite{GlM20}.}
Let $T:X \rightarrow X$ be a positive bounded linear operator on an ordered Banach space $(X,K)$. Further assume that $\inner\, K \neq \emptyset$. Then%
\begin{equation}\label{eq:Interior-estimate-Spectral-radius-new}
  r(T) \geq \inf\{ \lambda \in \R : \exists x \in \inner\, K \mbox{ s.t. } Tx \ll \lambda x \}.%
\end{equation}
\end{lemma}
%\NN{What does the symbol $\ll$ mean?}
\begin{IEEEproof}
Pick any $y \in \inner\,K$, any $\lambda>r(T)$, and define $x:= (\lambda I - T)^{-1}y$. Representing the resolvent by the Neumann series and using that $T$ is a positive operator, we obtain%
\begin{equation*}
  x = \sum_{k=0}^\infty \frac{T^k}{\lambda^{k+1}}y = \frac{1}{\lambda}y + \sum_{k=1}^\infty \frac{T^k}{\lambda^{k+1}}y \geq \frac{1}{\lambda}y,%
\end{equation*}
where we use that $r(T/\lambda) < 1$ which guarantees the convergence of the series. As $y \in \inner\,K$, this implies   $x \in \inner\,K$. Moreover,%
\begin{align*}
  Tx &= T(\lambda I - T)^{-1}y = (T - \lambda I + \lambda I)(\lambda I - T)^{-1}y\\
	   &= - y + \lambda(\lambda I - T)^{-1}y  = -y + \lambda x \ll \lambda x.%
\end{align*}
Hence, $(r(T),+\infty) \subset \{ \lambda \in \R : \exists x \in \inner(K) \mbox{ s.t. } Tx \leq \lambda x \}$, and the claim follows.
\end{IEEEproof}

\begin{IEEEbiography}[{\includegraphics[width=1in,height=1.3in,clip,keepaspectratio]{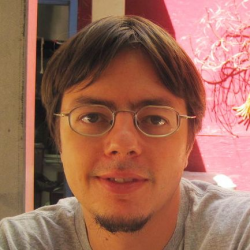}}]{Christoph Kawan} received the diploma and doctoral degree in mathematics (both under supervision of Prof. Fritz Colonius) 
from the University of Augsburg, Germany, in 2006 and 2009, respectively. He spent four months as a research scholar at the State University of Campinas, Brazil, in 2011, and nine months at the Courant Institute of Mathematical Sciences at New York University, in 2014. He is the author of the book `Invariance Entropy for Deterministic Control Systems - An Introduction' (Lecture Notes in Mathematics 2089. Berlin: Springer, 2013). He currently works as a researcher at the Ludwig-Maximilians-Universit\"{a}t in Munich. His research interests include networked and information-based control and non-autonomous dynamical systems.%
\end{IEEEbiography}

\begin{IEEEbiography}[{\includegraphics[width=1in,height=1.25in,clip,keepaspectratio]{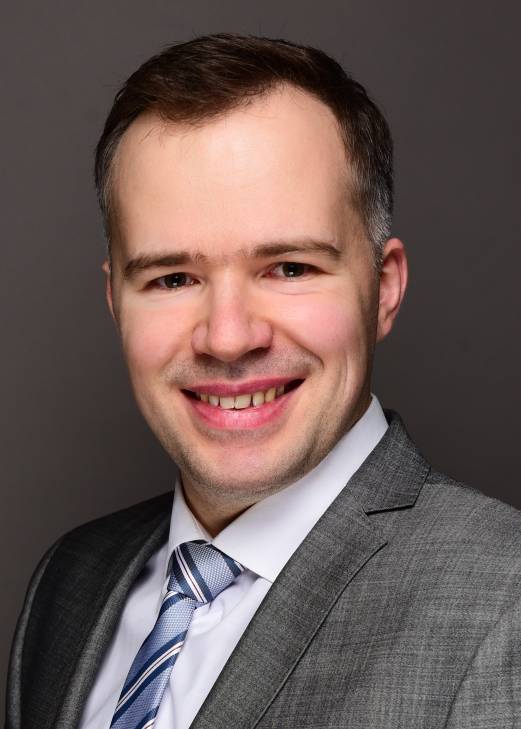}}]{Andrii Mironchenko}
was born in 1986. He received his MSc degree at the I.I. Mechnikov Odesa National University, Ukraine in 2008 and his PhD degree at the University of Bremen, Germany in 2012. 
He has held a research position at the University of W\"urzburg, Germany and was a Postdoctoral Fellow of Japan Society for Promotion of Science (JSPS) at the Kyushu Institute of Technology, Japan (2013--2014). In 2014 he joined the Faculty of Mathematics and Computer Science at the University of Passau, Germany. Dr. Mironchenko is a (co)-author of over 50 papers in journals and conferences in control theory and applied mathematics. His research interests include infinite-dimensional systems, stability theory, hybrid systems and applications of control theory to biological systems. 
\end{IEEEbiography}

\begin{IEEEbiography}[{\includegraphics[width=1in,height=1.25in,clip,keepaspectratio]{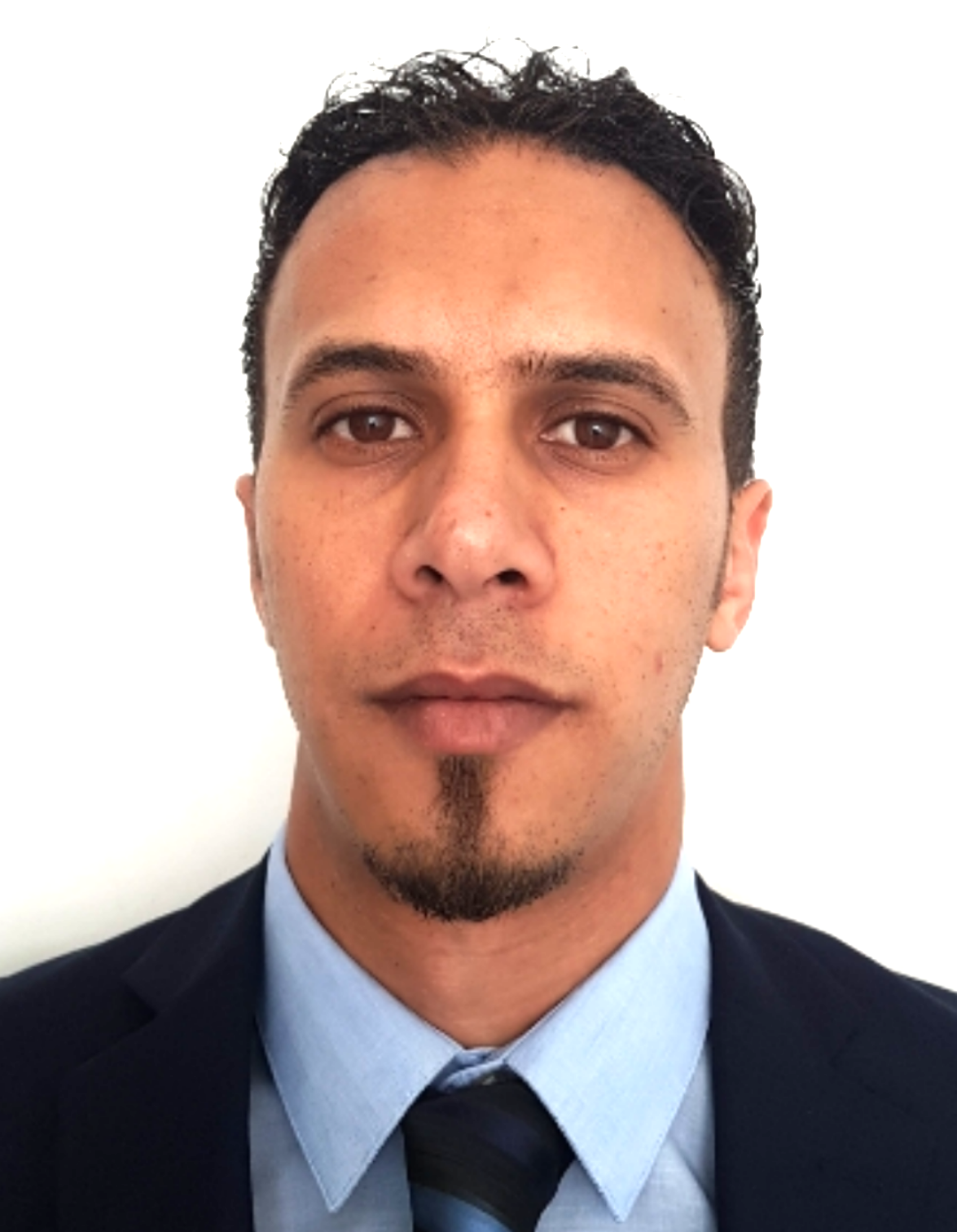}}]{Abdalla Swikir}
	is currently a PhD candidate in the Department of Electrical and Computer Engineering at Technical University of Munich (TUM), Germany, since November 2016.  He is also a faculty member in the Electrical Engineering Department, Omar Al-Mukhtar University, Albyada, Libya. He received a B.Sc. degree in Electrical Engineering in August 2008 from Omar Al-Mukhtar University, Albyada, Libya, and an M.Sc. degree in Electrical Engineering in December 2015 from Ohio State University, Columbus, OH, USA. From January 2010 to April 2013 he was a teaching assistant in the Department of Electrical Engineering at Omar Al-Mukhtar University, Albyada, Libya. His research interests are formal synthesis, symbolic models, and compositional methods of large-scale cyber-physical systems.
\end{IEEEbiography}

\begin{IEEEbiography}[{\includegraphics[width=1.1in,height=1.25in,clip,keepaspectratio]{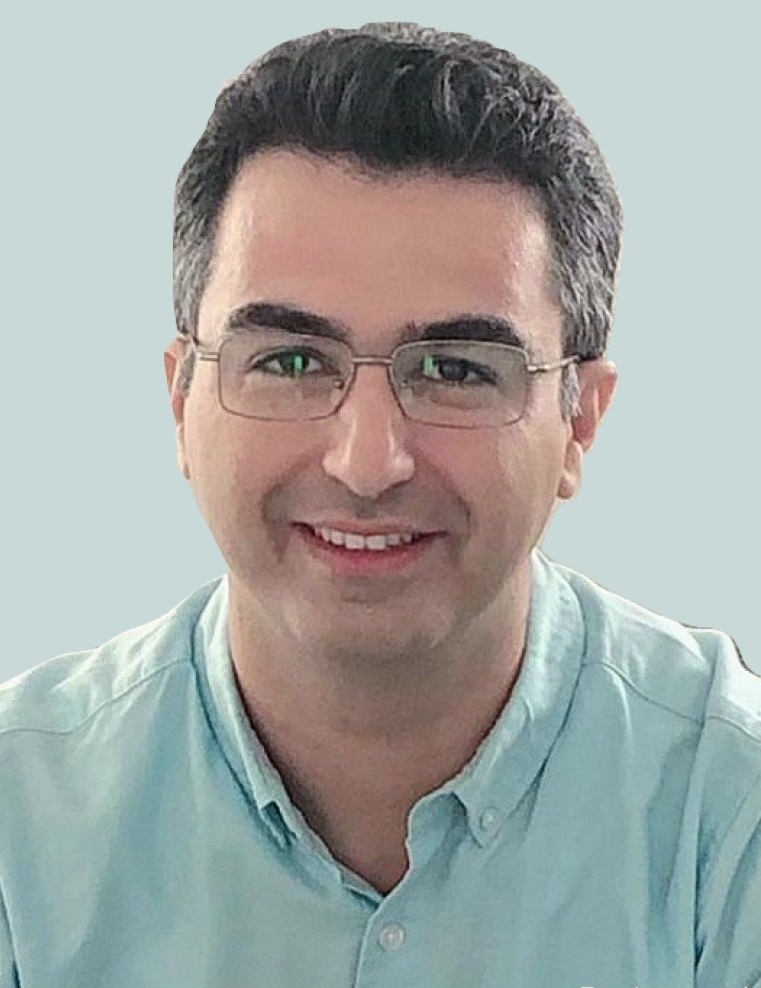}}]{Navid Noroozi}
received his M.Sc. and Ph.D. degrees in Electrical Engineering from Shiraz University, Shiraz, Iran, in 2009 and 2014, respectively.
From Aug 2012 to Feb 2013, he was a visiting scholar at the University of Melbourne, Melbourne, Australia. He joined Sheikh Bahaei University, Iran, in Aug 2014, where he was an Assistant Professor until Aug 2016. From Aug 2016 to June 2020 he has held postdoctoral positions at Passau University and Otto-von-Guericke University Magdeburg. Since July 2020 he has been a Scientific Associate at the Institute of Informatics, LMU Munich. Dr.~Noroozi was awarded a Humboldt Fellowship from the Alexander von Humboldt Foundation in 2015.
His research interests include large-scale complex networks, nonlinear systems and cyber-physical systems.
\end{IEEEbiography}

\begin{IEEEbiography}[{\includegraphics[width=1in,height=1.25in,clip,keepaspectratio]{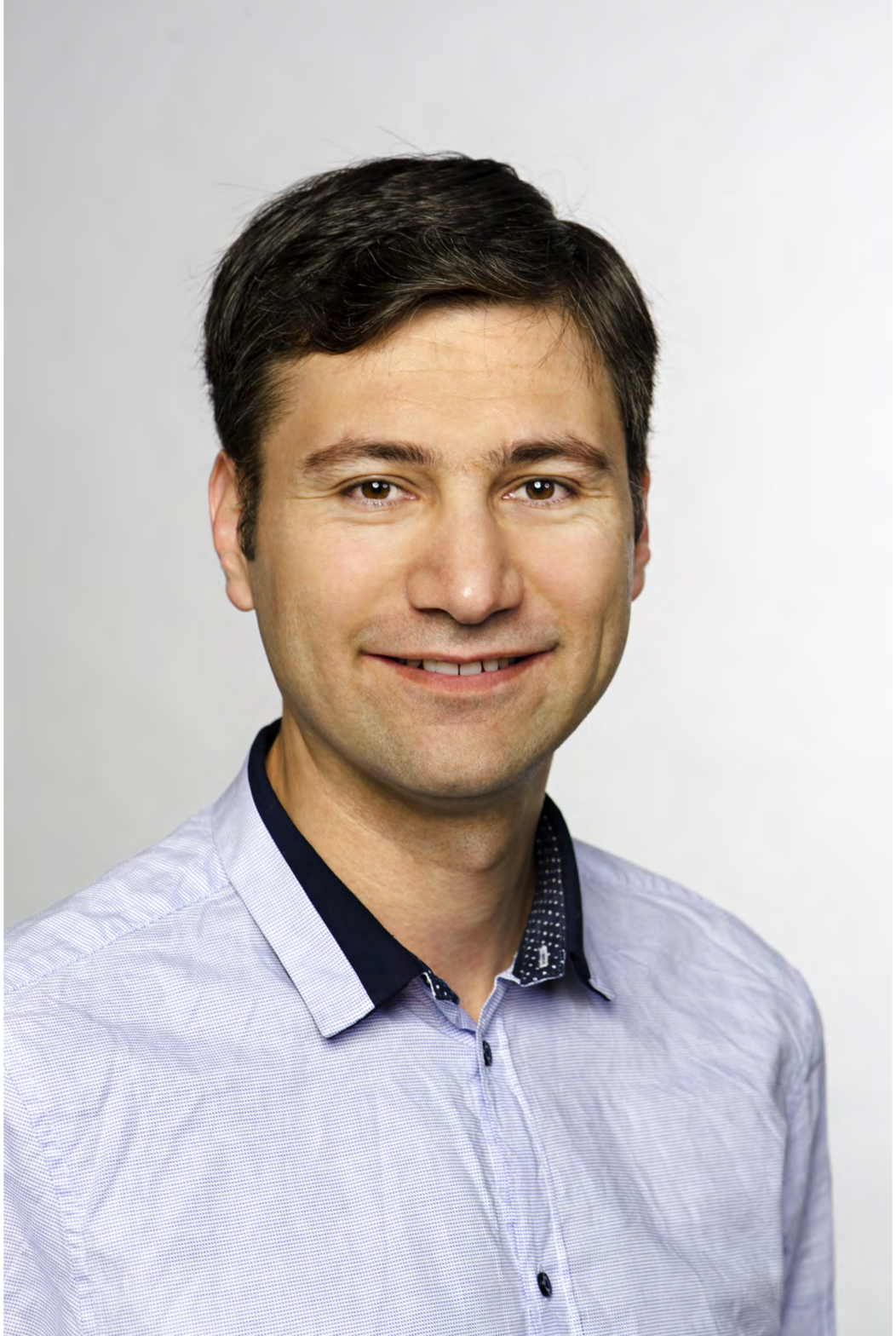}}]{Majid Zamani}
(M'12--SM'16)
is an Assistant Professor in the Computer Science Department at the University of Colorado Boulder, USA. He received a B.Sc. degree in Electrical Engineering in 2005 from Isfahan University of Technology, Iran, an M.Sc. degree in Electrical Engineering in 2007 from Sharif University of Technology, Iran, an MA degree in Mathematics and a Ph.D. degree in Electrical Engineering both in 2012 from University of California, Los Angeles, USA. Between September 2012 and December 2013 he was a postdoctoral researcher at the Delft Centre for Systems and Control, Delft University of Technology, Netherlands. From May 2014 to January 2019 he was an Assistant Professor in the Department of Electrical and Computer Engineering at the Technical University of Munich, Germany. From December 2013 to April 2014 he was an Assistant Professor in the Design Engineering Department, Delft University of Technology, Netherlands. He received an ERC starting grant award from the European Research Council in 2018.%

His research interests include verification and control of hybrid systems, embedded control software synthesis, networked control systems, and incremental properties of nonlinear control systems.
\end{IEEEbiography}

%>>>>>>> 5f968ccae93640ed0e9c350e78b160d149579522
\end{document}